\documentclass[12pt]{article}

\usepackage{amssymb}
\usepackage{amsfonts}
\usepackage{amsmath}
\usepackage{amsthm}
\usepackage[margin=0.70in]{geometry}
\usepackage{enumerate}
\usepackage{amstext}
\usepackage{layout}

\numberwithin{equation}{section}

\newtheorem{theorem}{Theorem}[section]

\newtheorem{lemma}{Lemma}[section]

\newtheorem{remark}{Remark}[section]

\begin{document}

\noindent
{\bf\large{Some results on probabilities of moderate deviations}}

\vskip 0.3cm

\noindent {\bf Deli Li\footnote{Deli Li, Department of Mathematical Sciences, Lakehead University, Thunder Bay, Ontario, Canada}
~~~Yu Miao\footnote{Yu Miao, College of Mathematics and Information Science,
Henan Normal University, Xinxiang, Henan, China}~~~
Yongcheng Qi\footnote{Yongcheng Qi,
Department of Mathematics and Statistics, University of Minnesota Duluth, Duluth, MN, USA}
}

\vskip 0.3cm

\noindent {\bf Abstract}~~Let $\{X, X_{n}; n \geq 1\}$ be a sequence of i.i.d. non-degenerate real-valued random variables
with $\mathbb{E}X^{2} < \infty$. Let $S_{n} = \sum_{i=1}^{n} X_{i}$, $n \geq 1$. Let
$g(\cdot): ~[0, \infty) \rightarrow [0, \infty)$ be a nondecreasing regularly varying function with index $\rho \geq 0$
and $\lim_{t \rightarrow \infty} g(t) = \infty$. Let $\mu = \mathbb{E}X$ and $\sigma^{2} = \mathbb{E}(X - \mu)^{2}$.
In this paper, on the scale $g(\log n)$, we obtain precise asymptotic estimates for the probabilities of
moderate deviations of the form $\displaystyle \log \mathbb{P}\left(S_{n} - n \mu > x \sqrt{ng(\log n)} \right)$,
$\displaystyle \log \mathbb{P}\left(S_{n} - n \mu < -x \sqrt{ng(\log n)} \right)$, and
$\displaystyle \log \mathbb{P}\left(\left|S_{n} - n \mu \right| > x \sqrt{ng(\log n)} \right)$ for all $x > 0$.
Unlike those known results in the literature, the moderate deviation results established in this paper depend
on both the variance and the asymptotic behavior of the tail distribution of $X$.

~\\


\vskip 0.3cm


\vskip 0.3cm

\noindent {\bf Running Head}:~~Probabilities of moderate deviations

\section{Introduction}

Throughout this paper, let $\{X, X_n; n \geq 1\}$ be a sequence of independent and identically distributed (i.i.d.)
real-valued random variables defined on a probability space $(\Omega , {\cal F}, \mathbb{P})$ and, as usual,
let $S_{n} = \sum_{i=1}^n X_i, n \geq 1$. Write $\log t = \log_{e} t$, $t > 0$ and define $\log 0 = - \infty$.

It is well known that Cram\'{e}r [10] and Chernoff [8] initiated the study of the theory of large deviations,
that characterizes the exponential concentration behaviour, as $n \rightarrow \infty$, of a sequence of
probabilities $\{\mathbb{P}(S_n / n \in \mathbf{A}) ; n \geq 1\}$, where $\mathbf{A} \subseteq (-\infty, \infty)$. They showed that
\begin{equation}
\left\{
\begin{array}{ll}
& \mbox{$\displaystyle \limsup_{n \to \infty} \frac{\log \mathbb{P}\left(S_{n}/n\in \mathbf{A} \right)}{n} \leq -\Lambda(\mathbf{A})~
\mbox{for every closed set}~\mathbf{A} \subseteq (-\infty, \infty),$}\\
&\\
& \mbox{$\displaystyle
\liminf_{n \to \infty} \frac{\log \mathbb{P}\left(S_{n}/n\in \mathbf{A} \right)}{n} \geq -\Lambda(\mathbf{A})~
\mbox{for every open set}~\mathbf{A} \subseteq (-\infty, \infty),$}
\end{array}
\right.
\end{equation}
provided that
\[
M(t) \equiv \mathbb{E}\left(e^{tX} \right) < \infty ~~\mbox{for all}~ t \in (-\infty, \infty),
\]
where, for $x \in (-\infty, \infty)$ and $\mathbf{A} \subseteq (-\infty, \infty)$,
\[
\Lambda(\mathbf{A}) = \inf_{x \in \mathbf{A}} I(x),~
I(x) = \sup_{t \in (-\infty, \infty)} \left(tx - \log M(t) \right).
\]
This fundamental result, which describes on the scale of a law of large number type ergodic phenomenon,
is what we call the Cram\'{e}r-Chernoff large deviation principle (in short, LDP) for $\{S_{n}; ~n \geq 1 \}$.
Clearly, the rate function $I(x), ~x \in (-\infty, \infty)$ of the LDP in (1.1) is determined by the moment generating function
$M(t), ~t \in (-\infty, \infty)$ of random variable $X$. The idea to study the log-moment generating function in a non i.i.d. setting
seems to go back  to Sievers [31], followed by Plachky [28] and Plachky and Steinebach [29]. More recently Comman [9] improved [29] substantially
by weakening the differentiability and convexity conditions. Donsker and Varadhan [13] and Bahadur and Zabell [1] established an LDP for
sums of i.i.d. Banach space-valued random variables. Bolthausen [3] extended the Cram\'{e}r-Chernoff-Donsker-Varadhan-Bahadur-Zabell
LDP when the laws of the random variables converge weakly and satisfy a uniform exponential integrability condition. As an application
of the Bolthausen LDP, Li, Rosalsky, and Al-Mutairi [22] established an LDP for bootstrapped sample means. Since large deviation theory
deals with the decay of the probability of increasingly unlikely events, it has applications in many different scientific fields, ranging
from queuing theory to statistics and from finance to engineering. There have been a great number of investigations on the probabilities
of large deviations for sums of independent random variables. Surveys of these investigations can be found in Book [4],
Dembo and Zeitouni [12], Petrov [26, 27]), Saulis and Statulevi\u{c}ius [30], Stroock [33], etc.

Inspired by the results Gantert [15] and Hu and Nyrhinen [16], Li and Miao [19] established an LDP, on the scale $\log n$,
for partial sums of i.i.d. $\mathbf{B}$-valued random variables. For the special case $\mathbf{B} = (- \infty, \infty)$, if
$S_n / n^{1/p} \rightarrow_{\mathbb{P}} 0$ for some $0<p<2$ then, for all $s > 0$,
\begin{equation}
\left\{
\begin{array}{ll}
& \mbox{$\displaystyle
\limsup_{n \to \infty} \frac{\log \mathbb{P}\left(\left|S_{n} \right| > s n^{1/p} \right)}{\log n}
= - (\overline{\beta} -p)/p,$}\\
&\\
& \mbox{$\displaystyle
\liminf_{n \to \infty} \frac{\log \mathbb{P}\left(\left|S_{n} \right| > s n^{1/p} \right)}{\log n}
= -( \underline{\beta} -p)/p,$}
\end{array}
\right.
\end{equation}
where
\[
\overline{\beta} = - \limsup_{t \rightarrow \infty}
\frac{\log \mathbb{P}(\log|X| > t)}{t}
~~\mbox{and}~~\underline{\beta} = - \liminf_{t \rightarrow \infty}
\frac{\log \mathbb{P}(\log|X| > t)}{t}.
\]
In particular, under the same hypotheses, for all $s>0$,
\[
\lim_{n \to \infty} \frac{\log \mathbb{P}\left(\left|S_{n} \right| > s n^{1/p} \right)}{\log n}
= - ({\beta} -p)/p ~~\mbox{if and only if}~~\overline{\beta}
= \underline{\beta}= \beta.
\]
As a special case of (1.2), the main results of Hu and Nyrhinen [16] are not only improved,
but also extended. Recently, a similar large deviation result to (1.2)
for partial sums of i.i.d. random variables with super-heavy tailed distribution is established
in Li, Miao, and Stoica [20] which extends in particular the results of Stoica [32] and Nakata [23].
We must point out that the large deviation results established in both Li and Miao [18] and Li, Miao, and Stoica [20]
depend only on the asymptotic behaviour of the tail distribution $\mathbb{P}(|X| > t)$ as $t \to \infty$.

In the another direction, under the Cram\'{e}r condition, which asserts that
\[
\mathbb{E}\left(e^{t|X|} \right) < \infty ~~\mbox{for some}~ t > 0,
\]
Petrov [24] obtained asymptotic expansions for
\[
\mathbb{P} \left(S_{n} > n \mu + n^{1/2}x \right) ~~\mbox{and}~~\mathbb{P} \left(S_{n} < n \mu - n^{1/2}x \right)
\]
for $x \geq 0$ and $x = {\it o}\left(n^{1/2} \right)$, where $\mu = \mathbb{E}X$. Let $\{b_{n};~n \geq 1 \}$ be a
sequence of positive real numbers such that
\[
\lim_{n \to \infty} \frac{b_{n}}{n} = 0~~\mbox{and}~~\lim_{n \to \infty} \frac{b_{n}}{\sqrt{n}} = \infty.
\]
Then it follows from the Petrov asymptotic expansions [24] that
\begin{equation}
\left\{
\begin{array}{ll}
& \mbox{$\displaystyle \limsup_{n \to \infty} \frac{n}{b_{n}^{2}} \log \mathbb{P}\left(\frac{S_{n} - n \mu}{b_{n}}
\in \mathbf{A} \right) \leq - \inf_{x \in \mathbf{A}}\frac{x^{2}}{2 \sigma^{2}}~
\mbox{for closed }~\mathbf{A} \subseteq (-\infty, \infty),$}\\
&\\
& \mbox{$\displaystyle \liminf_{n \to \infty} \frac{n}{b_{n}^{2}} \log \mathbb{P}\left(\frac{S_{n} - n \mu}{b_{n}}
\in \mathbf{A} \right) \geq - \inf_{x \in \mathbf{A}}\frac{x^{2}}{2 \sigma^{2}}~
\mbox{for open}~\mathbf{A} \subseteq (-\infty, \infty),$}\\
\end{array}
\right.
\end{equation}
where $\sigma^{2} = \mbox{Var}\left(X - \mathbb{E}X \right)^{2}$. This classical result, which describes the
probabilities on a scale between a law of large numbers and some sort of central limit theorem,
is what we call the moderate deviation principle (in short, MDP) for $\{S_{n}; ~n \geq 1 \}$. Under the assumptions that
\[
\mathbb{E}X = 0, ~\frac{b_{n}}{n} \downarrow 0,~~\mbox{and}~~\frac{b_{n}}{\sqrt{n}} \uparrow \infty ~~\mbox{as}~ n \to \infty,
\]
Eichelsbacher and L\"{o}we [14] showed that (1.3) holds for some $\sigma^{2} < \infty$ if and only if
\[
\sigma^{2} = \mathbb{E}X^{2}~~\mbox{and}~~\lim_{n \to \infty} \frac{n}{b_{n}^{2}} \log \left(n \mathbb{P}\left(|X| > b_{n} \right) \right) = - \infty.
\]
Clearly, the rate function $I(x) = \frac{x^{2}}{2 \sigma^{2}}, ~x \in (-\infty, \infty)$ of the MDP in (1.3) is determined
by the variance of random variable $X$. Borovkov and Mogul'ski\u{i} [5], Chen [6, 7], de Acosta [11], and Ledoux [17]
obtained versions of (1.3) in a Banach space setting under various conditions.

Motivated by (1.2), (1.3), and the work of Eichelsbacher and L\"{o}we [14], it is natural to ask if
any MDP result for partial sums of $\{X, X_{n};~n \geq 1 \}$ is determined by either the variance of $X$ only or
the tail distribution of $X$ only. In this paper, we focus on this problem to study the MDP under the condition
$\mathbb{E}X^{2} < \infty$ only. On the scale $g(\log n)$, we obtain precise asymptotic estimates for the probabilities of
moderate deviations of the form $\log \mathbb{P}\left(S_{n} - n \mu > x \sqrt{ng(\log n)} \right)$,
$\log \mathbb{P}\left(S_{n} - n \mu < -x \sqrt{ng(\log n)} \right)$, and
$\log \mathbb{P}\left(\left|S_{n} - n \mu \right| > x \sqrt{ng(\log n)} \right)$ for all $x > 0$, where
$g(\cdot)$: $[0, \infty) \rightarrow [0, \infty)$ is a non-decreasing regularly varying function with index $\rho \geq 0$
and $\lim_{t \rightarrow \infty} g(t) = \infty$. From the results established in this paper, we can see that, unlike
those known results, the moderate deviation results established in this paper depend on both the variance and the asymptotic
behavior of the tail distribution of the random variable $X$ and hence, the answer to the problem stated at the beginning of this
paragraph is negative.

The plan of the paper is as follows. Our main results Theorems 2.1 and 2.2, which are some general results on probabilities
of moderate deviations for partial sums of $\{X, X_{n};~n \geq 1 \}$, are presented in Section 2. Some preliminary results
needed to prove the main results are listed (and proved) in Section 3. The proofs of Theorems 2.1 and 2.2 are given in Sections 4
and 5 respectively. The truncation technique, conditional probability technique, two preliminary results on the regularly varying
functions, a maximal inequality, and Kolmogorov exponential inequalities are paramount in the proof of Theorem 2.1. The
main tools employed in proving Theorem 2.2 are the symmetrization technique, Kolmogorov strong law of large numbers,
Kolmogorov exponential inequalities, and the method of proof by contradiction.

\section{Statement of the main results}

Before we can formulate our results, we need some extra notation. For any real numbers $a$ and $b$, let
$a \wedge b = \min\{a, b\}$, $a \vee b = \max\{a, b\}$, $a^{+} = a \vee 0$, and
$a^{-} = (-a)\vee 0$.  For any real numbers $x \in (0, \infty)$ and $y \in (-\infty, \infty)$ put,
by convention, $(\pm\infty + y)/x =  \pm \infty$. Let $\rho \geq 0$ and let $\mathcal{V}_{\rho}$ be the set
of all nondecreasing regularly varying functions with index $\rho$ and $\lim_{t \to \infty} g(t) = \infty$.
Thus, if $g(\cdot) \in \mathcal{V}_{\rho}$, then
\[
\lim_{t \to \infty} \frac{g(xt)}{g(t)} = x^{\rho}~~\mbox{for all}~ x > 0.
\]
Clearly, if $g(\cdot) \in \mathcal{V}_{0}$, then $g(\cdot)$ is slowly varying at infinity; if $g(\cdot)$
is a nondecreasing regularly varying function with index $\rho > 0$, then $\lim_{t \to \infty} g(t) = \infty$
automatically. Let $X$ be a real-valued random variable. For any given $g(\cdot) \in \mathcal{V}_{\rho}$, write
\begin{equation}
\overline{\lambda}_{1}
= - \limsup_{t \rightarrow \infty} \frac{\log \left(t^{2} \mathbb{P}(X > t ) \right)}{g(\log t)}
~\mbox{and}~\underline{\lambda}_{1} = - \liminf_{t \rightarrow \infty}
\frac{\log \left(t^{2} \mathbb{P}(X > t ) \right)}{g(\log t)},
\end{equation}
\begin{equation}
\overline{\lambda}_{2}
= - \limsup_{t \rightarrow \infty} \frac{\log \left(t^{2} \mathbb{P}(X < -t ) \right)}{g(\log t)}
~\mbox{and}~\underline{\lambda}_{2} = - \liminf_{t \rightarrow \infty}
\frac{\log \left(t^{2} \mathbb{P}(X < -t ) \right)}{g(\log t)},
\end{equation}
and
\begin{equation}
\overline{\lambda}
= - \limsup_{t \rightarrow \infty} \frac{\log \left(t^{2} \mathbb{P}(|X| > t ) \right)}{g(\log t)}
~\mbox{and}~\underline{\lambda} = - \liminf_{t \rightarrow \infty}
\frac{\log \left(t^{2} \mathbb{P}(|X| > t ) \right)}{g(\log t)}.
\end{equation}
Clearly, $\overline{\lambda}_{1}$ and $\underline{\lambda}_{1}$ defined in (2.1) are two parameters of $X$
determined by the asymptotic behavior of the tail distribution $\mathbb{P}(X > t)$ as $t \to \infty$,
$\overline{\lambda}_{2}$ and $\underline{\lambda}_{2}$ defined in (2.2) are two parameters of $X$ determined by
the asymptotic behavior of the tail distribution $\mathbb{P}(X < -t)$ as $t \to \infty$, and $\overline{\lambda}$
and $\underline{\lambda}$ defined in (2.3) are two parameters of $X$ determined
by the asymptotic behavior of the tail distribution $\mathbb{P}(|X| > t)$ as $t \to \infty$.

The following Theorem 2.1 provides general and precise moderate deviation results for the partial sums of
$\{X, X_{n}; n \geq 1 \}$ under the finite second moment condition only.

\vskip 0.2cm

\begin{theorem}
Let $\{X, X_{n}; n \geq 1\}$ be a sequence of i.i.d. non-degenerate real-valued random variables with
$\mathbb{E}X^{2} < \infty$. Write $\mu = \mathbb{E}X$ and $\sigma^{2} = \mathbb{E}(X - \mu)^{2} \in (0, \infty)$.
Then, for any given $g(\cdot) \in \mathcal{V}_{\rho}$, we have
\begin{equation}
\left\{
\begin{array}{ll}
&
\mbox{$\displaystyle \limsup_{n \rightarrow \infty} \frac{\log \mathbb{P}\left(S_{n} -  n \mu > x \sqrt{n g(\log n)} \right)}{g(\log n)}
= - \left(\frac{x^{2}}{2\sigma^{2}} \wedge \frac{\overline{\lambda}_{1}}{2^{\rho}} \right)~\mbox{for all}~x > 0$,}\\
&\\
&
\mbox{$\displaystyle \liminf_{n \rightarrow \infty} \frac{\log \mathbb{P}\left(S_{n} - n \mu > x \sqrt{n g(\log n)} \right)}{g(\log n)}
= - \left(\frac{x^{2}}{2\sigma^{2}} \wedge \frac{\underline{\lambda}_{1}}{2^{\rho}} \right)~\mbox{for all}~x > 0$,}
\end{array}
\right.
\end{equation}
\begin{equation}
\left\{
\begin{array}{ll}
&
\mbox{$\displaystyle \limsup_{n \rightarrow \infty} \frac{\log \mathbb{P}\left(S_{n} - n \mu < -x \sqrt{n g(\log n)} \right)}{g(\log n)}
= - \left(\frac{x^{2}}{2\sigma^{2}} \wedge \frac{\overline{\lambda}_{2}}{2^{\rho}} \right)~\mbox{for all}~x > 0$,}\\
&\\
&
\mbox{$\displaystyle \liminf_{n \rightarrow \infty} \frac{\log \mathbb{P}\left(S_{n} - n \mu < - x \sqrt{n g(\log n)} \right)}{g(\log n)}
= - \left(\frac{x^{2}}{2\sigma^{2}} \wedge \frac{\underline{\lambda}_{2}}{2^{\rho}} \right)~\mbox{for all}~x > 0$,}
\end{array}
\right.
\end{equation}
and
\begin{equation}
\left\{
\begin{array}{ll}
&
\mbox{$\displaystyle
\limsup_{n \rightarrow \infty} \frac{\log \mathbb{P}\left(\left|S_{n} - n \mu \right| > x \sqrt{n g(\log n)} \right)}{g(\log n)}
= - \left(\frac{x^{2}}{2\sigma^{2}} \wedge \frac{\overline{\lambda}}{2^{\rho}} \right)~\mbox{for all}~x > 0$,}\\
&\\
&
\mbox{$\displaystyle
\liminf_{n \rightarrow \infty} \frac{\log \mathbb{P}\left(\left|S_{n} - n \mu \right| > x \sqrt{n g(\log n)} \right)}{g(\log n)}
= - \left(\frac{x^{2}}{2\sigma^{2}} \wedge \frac{\underline{\lambda}}{2^{\rho}} \right)~\mbox{for all}~x > 0$.}
\end{array}
\right.
\end{equation}
Hence,
\[
\lim_{n \rightarrow \infty} \frac{\log \mathbb{P}\left(S_{n}- n \mu > x \sqrt{n g(\log n)} \right)}{g(\log n)}
= - \left(\frac{x^{2}}{2\sigma^{2}} \wedge \frac{\hat{\lambda}_{1}}{2^{\rho}} \right)~~\mbox{for all}~ x > 0~
\mbox{if and only if}~ \overline{\lambda}_{1} = \underline{\lambda}_{1} = \hat{\lambda}_{1},
\]
\[
\lim_{n \rightarrow \infty} \frac{\log \mathbb{P}\left(S_{n} - n \mu < - x \sqrt{n g(\log n)} \right)}{g(\log n)}
= - \left(\frac{x^{2}}{2\sigma^{2}} \wedge \frac{\hat{\lambda}_{1}}{2^{\rho}} \right)~~\mbox{for all}~ x > 0~
\mbox{if and only if}~ \overline{\lambda}_{2} = \underline{\lambda}_{2} = \hat{\lambda}_{2},
\]
and
\[
\lim_{n \rightarrow \infty} \frac{\log \mathbb{P}\left(\left|S_{n} - n \mu \right| > x \sqrt{n g(\log n)} \right)}{g(\log n)}
= - \left(\frac{x^{2}}{2\sigma^{2}} \wedge \frac{\hat{\lambda}}{2^{\rho}} \right)~~\mbox{for all}~ x > 0~
\mbox{if and only if}~ \overline{\lambda} = \underline{\lambda} = \hat{\lambda}.
\]
In particular,
\[
\lim_{n \rightarrow \infty} \frac{\log \mathbb{P}\left(S_{n} - n \mu > x \sqrt{n g(\log n)} \right)}{g(\log n)}
= - \frac{x^{2}}{2\sigma^{2}}~~\mbox{for all}~ x > 0~
\mbox{if and only if}~ \overline{\lambda}_{1} = \infty,
\]
\[
\lim_{n \rightarrow \infty} \frac{\log \mathbb{P}\left(S_{n} - n \mu < - x \sqrt{n g(\log n)} \right)}{g(\log n)}
= - \frac{x^{2}}{2\sigma^{2}}~~\mbox{for all}~ x > 0~
\mbox{if and only if}~ \overline{\lambda}_{2} = \infty,
\]
and
\[
\lim_{n \rightarrow \infty} \frac{\log \mathbb{P}\left(\left|S_{n} - n \mu \right| > x \sqrt{n g(\log n)} \right)}{g(\log n)}
= - \frac{x^{2}}{2\sigma^{2}}~~\mbox{for all}~ x > 0~
\mbox{if and only if}~ \overline{\lambda} = \infty.
\]
\end{theorem}

\vskip 0.2cm

\begin{remark}
{\bf (i)}~~It is interesting to see that, on the scale $g(\log n)$, Theorem 2.1 provides us precise
asymptotic estimates for the probabilities of moderate deviations of
$\log \mathbb{P}\left(S_{n} - n \mu > x \sqrt{ng(\log n)} \right)$,
$\log \mathbb{P}\left(S_{n} - n \mu < -x \sqrt{ng(\log n)} \right)$, and
$\log \mathbb{P}\left(\left|S_{n} - n \mu \right| > x \sqrt{ng(\log n)} \right)$ for all $x > 0$
and such moderate deviation results depend only on both the variance and the asymptotic behavior
of the tail distribution of $X$.

\vskip 0.2cm

\noindent

{\bf (ii)}~~Since, for all $t > 0$,
\[
\begin{array}{lll}
\mbox{$\displaystyle
\frac{\log \left(t^{2} \mathbb{P}(X > t ) \right)}{g(\log t)} \vee \frac{\log \left(t^{2} \mathbb{P}(X < -t ) \right)}{g(\log t)}$}
 & \leq &
 \mbox{$\displaystyle
\frac{\log \left(t^{2} \mathbb{P}(|X| > t ) \right)}{g(\log t)}$}\\
&&\\
& \leq &
\mbox{$\displaystyle
\frac{\log \left(2 \left(t^{2} \mathbb{P}(X > t ) \vee t^{2} \mathbb{P}(X < -t ) \right)
\right)}{g(\log t)}$}\\
&&\\
& = &
\mbox{$\displaystyle
\frac{\log 2}{g(\log t)} + \left(
\frac{\log \left(t^{2} \mathbb{P}(X > t ) \right)}{g(\log t)} \vee \frac{\log \left(t^{2} \mathbb{P}(X < -t ) \right)}{g(\log t)}
\right),$}
\end{array}
\]
we have
\[
-\overline{\lambda}  = \left(-\overline{\lambda}_{1}\right) \vee \left(-\overline{\lambda}_{2} \right);
~\mbox{i.e.,}~~\overline{\lambda}  = \overline{\lambda}_{1} \wedge \overline{\lambda}_{2}.
\]
Similarly, we can show that
\[
\left(-\underline{\lambda}_{1}\right) \vee \left(-\underline{\lambda}_{2} \right) \leq -\underline{\lambda};
~\mbox{i.e.,}~~\underline{\lambda}  \leq \underline{\lambda}_{1} \wedge \underline{\lambda}_{2}.
\]
However, the assertion $\underline{\lambda} = \underline{\lambda}_{1} \wedge \underline{\lambda}_{2} $
is not true.
\end{remark}

\vskip 0.2cm

\begin{remark}
If, for the given $g(\cdot) \in \mathcal{V}_{\rho}$, $\lim_{t \to \infty} \frac{g(t)}{t} = \infty$ (for such case, $\rho \geq 1$),
then $\lim_{t \to \infty} \frac{\log t^{2}}{g(\log t)} = 0$ and hence, one can easily see that the $\overline{\lambda}_{1}$,
$\underline{\lambda}_{1}$, $\overline{\lambda}_{2}$, $\underline{\lambda}_{2}$, $\overline{\lambda}$,
and $\underline{\lambda}$ can also be defined respectively by
\[
\overline{\lambda}_{1}
= - \limsup_{t \rightarrow \infty} \frac{\log \mathbb{P}(X > t )}{g(\log t)},
~~\underline{\lambda}_{1} = - \liminf_{t \rightarrow \infty}
\frac{\log \mathbb{P}(X > t )}{g(\log t)},
\]
\[
\overline{\lambda}_{2}
= - \limsup_{t \rightarrow \infty} \frac{\log \mathbb{P}(X < -t )}{g(\log t)},
~~\underline{\lambda}_{2} = - \liminf_{t \rightarrow \infty}
\frac{\log \mathbb{P}(X < -t )}{g(\log t)},
\]
\[
\overline{\lambda}
= - \limsup_{t \rightarrow \infty} \frac{\log \mathbb{P}(|X| > t )}{g(\log t)},
~\mbox{and}~~\underline{\lambda} = - \liminf_{t \rightarrow \infty}
\frac{\log \mathbb{P}(|X| > t )}{g(\log t)}.
\]
\end{remark}

\vskip 0.2cm

\begin{remark}
Under the assumptions of Theorem 2.1, one can show that the $\overline{\lambda}_{1}$,
$\underline{\lambda}_{1}$, $\overline{\lambda}_{2}$, $\underline{\lambda}_{2}$, $\overline{\lambda}$,
and $\underline{\lambda}$ can also be defined respectively by
\[
\overline{\lambda}_{1} = \sup\left \{r \geq 0: \lim_{t \to \infty} t^{2}e^{r g(\log t)}\mathbb{P}(X > t) = 0 \right \},
~\underline{\lambda}_{1} = \sup\left \{r \geq 0: \liminf_{t \to \infty} t^{2}e^{r g(\log t)}\mathbb{P}(X > t) = 0 \right \},
\]
\[
\overline{\lambda}_{2} = \sup\left \{r \geq 0: \lim_{t \to \infty} t^{2}e^{r g(\log t)}\mathbb{P}(X < -t) = 0 \right \},
~\underline{\lambda}_{2} = \sup\left \{r \geq 0: \liminf_{t \to \infty} t^{2}e^{r g(\log t)}\mathbb{P}(X < -t) = 0 \right \},
\]
\[
\overline{\lambda} = \sup\left \{r \geq 0: \lim_{t \to \infty} t^{2}e^{r g(\log t)}\mathbb{P}(|X| > t) = 0 \right \},
~\mbox{and}~ \underline{\lambda} = \sup\left \{r \geq 0: \liminf_{t \to \infty} t^{2}e^{r g(\log t)}\mathbb{P}(|X| > t) = 0 \right \}.
\]
The proofs are left to the reader.
\end{remark}

\vskip 0.2cm

\begin{remark}
Let $\{X, X_{n}; n \geq 1\}$ be a sequence of i.i.d. non-degenerate real-valued random variables with $\mathbb{E}X^{2} < \infty$.
Then, following from Theorem 2.1 and Remark 2.3, we have the following two special cases which are related to the law of the iterated
logarithm for partial sums of i.i.d. random variables.

\vskip 0.2cm

\noindent
{\bf (i)}~~If $g(t) = t$, $t \geq 0$, then all conclusions in Theorem 2.1 hold with $g(\log n) = \log n$, $\rho = 1$, and
$\overline{\lambda}_{1}$, $\underline{\lambda}_{1}$, $\overline{\lambda}_{2}$, $\underline{\lambda}_{2}$,
$\overline{\lambda}$, and $\underline{\lambda}$ defined by
\[
\overline{\lambda}_{1} = \sup\left \{r \geq 0: \lim_{t \to \infty} t^{2+r}\mathbb{P}(X > t) = 0 \right \},
~\underline{\lambda}_{1} = \sup\left \{r \geq 0: \liminf_{t \to \infty} t^{2+r}\mathbb{P}(X > t) = 0 \right \},
\]
\[
\overline{\lambda}_{2} = \sup\left \{r \geq 0: \lim_{t \to \infty} t^{2+r}\mathbb{P}(X < -t) = 0 \right \},
~\underline{\lambda}_{2} = \sup\left \{r \geq 0: \liminf_{t \to \infty} t^{2+r}\mathbb{P}(X < -t) = 0 \right \},
\]
\[
\overline{\lambda} = \sup\left \{r \geq 0: \lim_{t \to \infty} t^{2+r}\mathbb{P}(|X| > t) = 0 \right \},
~\mbox{and}~ \underline{\lambda} = \sup\left \{r \geq 0: \liminf_{t \to \infty} t^{2+2}\mathbb{P}(|X| > t) = 0 \right \}.
\]
In particular,
\[
\lim_{n \rightarrow \infty} \frac{\log \mathbb{P}\left(S_{n} > x \sqrt{n \log n} \right)}{\log n}
= - \frac{x^{2}}{2\sigma^{2}}~~\mbox{for all}~ x > 0~
\mbox{if and only if}~ \mathbb{E}\left(X^{+} \right)^{r} < \infty~~\mbox{for all}~ r > 0,
\]
\[
\lim_{n \rightarrow \infty} \frac{\log \mathbb{P}\left(S_{n} < - x \sqrt{n \log n} \right)}{\log n}
= - \frac{x^{2}}{2\sigma^{2}}~~\mbox{for all}~ x > 0~
\mbox{if and only if}~ \mathbb{E}\left(X^{-} \right)^{r} < \infty~~\mbox{for all}~ r > 0,
\]
and
\[
\lim_{n \rightarrow \infty} \frac{\log \mathbb{P}\left(\left|S_{n} \right| > x \sqrt{n \log n} \right)}{\log n}
= - \frac{x^{2}}{2\sigma^{2}}~~\mbox{for all}~ x > 0~
\mbox{if and only if}~ \mathbb{E}|X|^{r} < \infty~~\mbox{for all}~ r > 0.
\]

\vskip 0.2cm

\noindent
{\bf (ii)}~~If $g(t) = \log (t \vee 1)$, $t \geq 0$, then all conclusions in Theorem 2.1 hold
with $g(\log n) = \log \log n$ ($n \geq 3$), $\rho = 0$, and $\overline{\lambda}_{1}$, $\underline{\lambda}_{1}$,
$\overline{\lambda}_{2}$, $\underline{\lambda}_{2}$, $\overline{\lambda}$, and $\underline{\lambda}$ defined by
\[
\overline{\lambda}_{1} = \sup\left \{r \geq 0: \lim_{t \to \infty} t^{2}(\log t)^{r} \mathbb{P}(X > t) = 0 \right \},
~\underline{\lambda}_{1} = \sup\left \{r \geq 0: \liminf_{t \to \infty} t^{2}(\log t)^{r}\mathbb{P}(X > t) = 0 \right \},
\]
\[
\overline{\lambda}_{2} = \sup\left \{r \geq 0: \lim_{t \to \infty}t^{2}(\log t)^{r}\mathbb{P}(X < -t) = 0 \right \},
~\underline{\lambda}_{2} = \sup\left \{r \geq 0: \liminf_{t \to \infty}t^{2}(\log t)^{r} \mathbb{P}(X < -t) = 0 \right \},
\]
\[
\overline{\lambda} = \sup\left \{r \geq 0: \lim_{t \to \infty} t^{2}(\log t)^{r} \mathbb{P}(|X| > t) = 0 \right \},
~\mbox{and}~ \underline{\lambda} = \sup\left \{r \geq 0: \liminf_{t \to \infty} t^{2}(\log t)^{r} \mathbb{P}(|X| > t) = 0 \right \}.
\]
In particular
\[
\lim_{n \rightarrow \infty} \frac{\log \mathbb{P}\left(S_{n} > x \sqrt{n \log \log n} \right)}{\log \log n}
= - \frac{x^{2}}{2\sigma^{2}}~~\mbox{for all}~ x > 0~
\mbox{if and only if}~ \mathbb{E}\left((X^{+})^{2} (\log X^{+})^{r}  \right) < \infty~~\mbox{for all}~ r > 0,
\]
\[
\lim_{n \rightarrow \infty} \frac{\log \mathbb{P}\left(S_{n} < -x \sqrt{n \log \log n} \right)}{\log \log n}
= - \frac{x^{2}}{2\sigma^{2}}~~\mbox{for all}~ x > 0~
\mbox{if and only if}~ \mathbb{E}\left((X^{-})^{2} (\log X^{-})^{r}  \right) < \infty~~\mbox{for all}~ r > 0,
\]
and
\[
\lim_{n \rightarrow \infty} \frac{\log \mathbb{P}\left(\left|S_{n}\right| > x \sqrt{n \log \log n} \right)}{\log \log n}
= - \frac{x^{2}}{2\sigma^{2}}~~\mbox{for all}~ x > 0~
\mbox{if and only if}~ \mathbb{E}\left(X^{2} (\log |X|)^{r}  \right) < \infty~~\mbox{for all}~ r > 0.
\]
\end{remark}

\vskip 0.2cm

The following Theorem 2.2 shows that $0 < \sigma^{2} = \mathbb{E}(X - \mu)^{2} < \infty$
is necessary for the moderate deviation results established in Theorem 2.1.

\vskip 0.2cm

\begin{theorem}
Let $\{X, X_{n}; n \geq 1\}$ be a sequence of i.i.d. real-valued random variables. Then, for any
given $\eta \in (-\infty, \infty)$ and $g(\cdot) \in \mathcal{V}_{\rho}$ for some $\rho \geq 0$, we have:

\vskip 0.2cm

\noindent
{\bf (i)}~~The following three statements are equivalent:
\begin{equation}
\lim_{n \rightarrow \infty} \frac{\log \mathbb{P}\left(\left|S_{n} - n \eta \right| > x \sqrt{n g(\log n)} \right)}{g(\log n)} = 0
~~\mbox{for all} ~x > 0;
\end{equation}
\begin{equation}
\lim_{n \rightarrow \infty} \frac{\log \mathbb{P}\left(\left|S_{n} - n \eta \right| > x \sqrt{n g(\log n)} \right)}{g(\log n)} = 0
~~\mbox{for some} ~x > 0;
\end{equation}
\begin{equation}
\mbox{Either} ~~\mathbb{E}X \neq \eta ~~\mbox{or}~~\mathbb{E}X^{2} = \infty~~\mbox{or}~ \sigma^{2} = \mbox{Var}(X) \in (0, \infty) ~\mbox{and}~~\lambda_{1} = \lambda_{2} = 0.
\end{equation}

\vskip 0.2cm

\noindent
{\bf (ii)}~~The following three statements are equivalent:
\begin{equation}
-\infty < \limsup_{n \rightarrow \infty} \frac{\log \mathbb{P}\left(\left|S_{n} - n \eta \right| > x \sqrt{n g(\log n)} \right)}{g(\log n)} < 0
~~\mbox{for all} ~x > 0;
\end{equation}
\begin{equation}
-\infty < \limsup_{n \rightarrow \infty} \frac{\log \mathbb{P}\left(\left|S_{n} - n \eta \right| > x \sqrt{n g(\log n)} \right)}{g(\log n)} < 0
~~\mbox{for some} ~x > 0;
\end{equation}
\begin{equation}
\sigma^{2} = \mbox{Var}(X) \in (0, \infty)~ \mbox{and}~\lambda_{1} > 0.
\end{equation}

\vskip 0.2cm

\noindent
{\bf (iii)}~~The following three statements are equivalent:
\begin{equation}
-\infty < \liminf_{n \rightarrow \infty} \frac{\log \mathbb{P}\left(\left|S_{n} - n \eta \right| > x \sqrt{n g(\log n)} \right)}{g(\log n)} < 0
~~\mbox{for all} ~x > 0;
\end{equation}
\begin{equation}
-\infty < \liminf_{n \rightarrow \infty} \frac{\log \mathbb{P}\left(\left|S_{n} - n \eta \right| > x \sqrt{n g(\log n)} \right)}{g(\log n)} < 0
~~\mbox{for some} ~x > 0;
\end{equation}
\begin{equation}
\sigma^{2} = \mbox{Var}(X) \in (0, \infty)~ \mbox{and}~\lambda_{2} > 0.
\end{equation}
\end{theorem}

\section{Preliminary lemmas}

In this section, we collect four preliminary lemmas needed for the proofs of our main results.
We need some additional notation. Let $m(Y)$ denote a median for a real-valued random variable $Y$.
We put $m(-Y) = -m(Y)$.

The following lemma is used to prove Theorem 2.1. The first part was obtained by Li and Rosalsky [20, Lemma 3.1]
and the second part was established by Petrov [25] (also see Petrov [27, Theorem 2.1]).

\vskip 0.2cm

\begin{lemma}
    Let $\{V_{k};~1 \leq k \leq n \}$ be a finite sequence of independent real-valued random variables and set
    $T_{0} = 0$ and $T_{k} = V_{1} + \cdots + V_{k}$, $1 \leq k \leq n$. Then, for every real $t$,
    \[
    \mathbb{P}\left(\max_{1 \leq k \leq n} \left(V_{k} + m\left(T_{k-1}\right) \right) > t \right)
    \leq 2 \mathbb{P}\left(\max_{1 \leq k \leq n} T_{k} > t \right),
    \]
    \[
    \mathbb{P}\left(\max_{1 \leq k \leq n} \left(T_{k} + m\left(T_{n} - T_{k}\right) \right) > t \right)
    \leq 2 \mathbb{P}\left(T_{n} > t \right).
    \]
\end{lemma}

\vskip 0.2cm

 \begin{lemma}
 Let $Y$ be a non-negative random variable. Let $p(\cdot): ~[t_{1}, \infty) \rightarrow [0, \infty)$
 be a non-decreasing function such that
 \begin{equation}
 0 < p(2t) \leq b p(t), ~t \geq t_{1}~~\mbox{and}~~\lim_{t \to \infty} p(t) \mathbb{P}(Y > t) = 0,
 \end{equation}
where $b > 1$ and $t_{1} > 0$ are two constants. Let $h(\cdot): ~[t_{2}, \infty) \rightarrow [0, \infty)$
be a non-decreasing function such that
 \begin{equation}
 \lim_{t \to \infty} h(t) = \infty ~~\mbox{and}~~\lim_{t \to \infty} \frac{h(t+1)}{h(t)} = 1,
 \end{equation}
where $t_{2} > 0$ is a constant. Then
 \begin{equation}
 \limsup_{n \to \infty} \frac{\log \left(p(n) \mathbb{P}\left(Y > n \right)\right)}{h(n)} = \limsup_{t \to \infty}
 \frac{\log \left(p(t) \mathbb{P}\left(Y > t \right)\right)}{h(t)}
 \end{equation}
 and
 \begin{equation}
 \liminf_{n \to \infty}
 \frac{\log \left(p(n) \mathbb{P}\left(Y > n \right)\right)}{h(n)} = \liminf_{t \to \infty}
 \frac{\log \left(p(t) \mathbb{P}\left(Y > t \right)\right)}{h(t)}.
 \end{equation}
 \end{lemma}

 \vskip 0.2cm

 \noindent {\it Proof}~~Since $p(\cdot): ~[t_{1}, \infty) \rightarrow [0, \infty)$
 is a non-decreasing function with (3.1), we have,
for $n \leq t < n + 1$ and all sufficiently large $n$,
 \[
 \begin{array}{lll}
 \mbox{$\displaystyle
\frac{1}{b} p(n+1)\mathbb{P}(Y > n+1)$}
& \leq & \mbox{$\displaystyle p(n) \mathbb{P}(Y > n+1)$}\\
&&\\
& \leq & \mbox{$\displaystyle
p(t) \mathbb{P}(Y > t) \leq p(n+1) \mathbb{P}(Y > n)$}\\
&&\\
& \leq & \mbox{$\displaystyle
bp(n) \mathbb{P}(Y > n)$}\\
&&\\
& \leq & 1.
\end{array}
 \]
Hence, for $n \leq t < n+1$ and all sufficiently large $n$,
 \[
 \begin{array}{lll}
 0 & \leq &
 \mbox{$\displaystyle  - \log \left(b p(n) \mathbb{P}(Y > n)\right)$}\\
 &&\\
 & \leq & \mbox{$\displaystyle - \log \left(p(t) \mathbb{P}(Y > t) \right)$}\\
  &&\\
& \leq & \mbox{$\displaystyle - \log \left(\frac{1}{b} p(n+1)\mathbb{P}(Y > n + 1)\right).$}
\end{array}
 \]
Since $h(\cdot): ~[t_{2}, \infty) \rightarrow [0, \infty)$ is a non-decreasing function,
we have, for $n \leq t < n + 1$ and all sufficiently large $n$,
\begin{equation}
\begin{array}{lll}
0 & \leq & \mbox{$\displaystyle
\frac{h(n)}{h(n+1)} \left(- \frac{\log b + \log \left(p(n) \mathbb{P}(Y > n)\right)}{h(n)}\right)$}\\
&&\\
& \leq & \mbox{$\displaystyle
- \frac{\log \left(p(t) \mathbb{P}(Y > t)\right)}{h(t)} $}\\
&&\\
& \leq & \mbox{$\displaystyle
 \frac{h(n+1)}{h(n)} \left(- \frac{- \log b + \log \left( p(n+1) \mathbb{P}(Y > n+1)\right)}{h(n+1)}\right).$}
 \end{array}
\end{equation}
Thus, (3.3) and (3.4) follow from (3.5) and (3.2). ~$\Box$

\vskip 0.2cm

\begin{lemma}
Let $X$ be a real-valued random variable with $\mathbb{E}X^{2} < \infty$. Then,
for any given $g(\cdot) \in \mathcal{V}_{\rho}$ and all $s > 0$, we have
\begin{equation}
\left\{
\begin{array}{ll}
& \mbox{$\displaystyle
\limsup_{n \to \infty} \frac{\log \left(n \mathbb{P}\left(X > s \sqrt{n g(\log n)} \right) \right)}{g(\log n)}
= - \frac{\underline{\lambda}_{1}}{2^{\rho}},$}\\
&\\
& \mbox{$\displaystyle
~\liminf_{n \to \infty} \frac{\log \left(n \mathbb{P}\left(X > s \sqrt{n g(\log n)} \right) \right)}{g(\log n)}
= - \frac{\underline{\lambda}_{1}}{2^{\rho}}$}
\end{array}
\right.
\end{equation}
and
\begin{equation}
\left\{
\begin{array}{ll}
& \mbox{$\displaystyle
\limsup_{n \to \infty} \frac{\log \left(n \mathbb{P}\left(X > \frac{s \sqrt{n}}{g(\log n)} \right) \right)}{g(\log n)}
= - \frac{\overline{\lambda}_{1}}{2^{\rho}},$}\\
&\\
& \mbox{$\displaystyle
~\liminf_{n \to \infty} \frac{\log \left(n \mathbb{P}\left(X > \frac{s \sqrt{n}}{g(\log n)} \right) \right)}{g(\log n)}
= - \frac{\underline{\lambda}_{1}}{2^{\rho}},$}
\end{array}
\right.
\end{equation}
where $\overline{\lambda}_{1} $ and $\underline{\lambda}_{1}$ are defined by (2.1).
\end{lemma}

\vskip 0.2cm

\noindent {\it Proof}~~Since $g(\cdot) \in \mathcal{V}_{\rho}$, we see that
\[
\tilde{g}(t) = g(n-1) + (t - n + 1)\left(g(n) - g(n-1) \right), ~n-1 \leq t < n,~ n \geq 1.
\]
is a continuous and nondecreasing function defined on $[0, \infty)$ such that
\[
\tilde{g}(n) = g(n), ~n \geq 1 ~\mbox{and}~~\lim_{t \to \infty} \frac{\tilde{g}(t)}{g(t)} = 1
\]
and hence, $\tilde{g}(\cdot) \in \mathcal{V}_{\rho}$. Thus, without loss of generality, we can assume
that $g(\cdot): [0, \infty) \to [0, \infty)$ is continuous (otherwise, $g(\cdot)$
can be replaced by $\tilde{g}(\cdot)$) and $g(1) > 0$ (since $\lim_{t \to \infty} g(t) = \infty$).

We first establish (3.6). For given $s > 0$, write
\[
\varphi_{s}(t) = s \sqrt{t g(\log (t \vee e))}, ~t \geq 0.
\]
Then, under the given conditions of $g(\cdot)$, $\varphi_{s}(\cdot):~[0, \infty) \rightarrow [0, \infty)$
is a continuous and strictly increasing function with $\varphi_{s}(0) = 0$ and $\lim_{t \to \infty} \varphi_{s}(t) = \infty$
and
\begin{equation}
\lim_{t \to \infty} g(\log t) = \infty, ~\lim_{t \to \infty} \frac{g(\log (t+1))}{g(\log t)} =1,
\end{equation}
\begin{equation}
\lim_{t \to \infty} \frac{g\left(\log \varphi_{s}(t)\right)}{g(\log t)}
= \lim_{t \to \infty} \frac{g\left(\frac{1}{2} \log t + \log s + \frac{1}{2} \log g(\log t) \right)}{g(\log t)}
= \frac{1}{2^{\rho}}.
\end{equation}
Let $\varphi_{s}^{-1}(\cdot)$ be the inverse function of $\varphi_{s}(\cdot)$. Under the given conditions, by Markov's inequality, we have
\[
\begin{array}{lll}
\mbox{$\displaystyle
\limsup_{t \to \infty} t \mathbb{P}\left(\varphi_{s}^{-1}\left(X^{+}\right) > t \right)$}
& = &
\mbox{$\displaystyle
\limsup_{t \to \infty} t \mathbb{P}\left(X > \varphi_{s}(t) \right) $}\\
&&\\
& \leq &
\mbox{$\displaystyle
\limsup_{t \to \infty} \frac{t \mathbb{E}X^{2}}{s^{2}tg(\log(t\vee e))}$}\\
&&\\
& = &
\left(\frac{\mathbb{E}X^{2}}{s^{2}} \right) \mbox{$\displaystyle
	\limsup_{t \to \infty} \frac{1}{g(\log(t\vee e))}$}\\
&&\\
& = & 0.
\end{array}
\]
Thus
\begin{equation}
\lim_{t \to \infty} t \mathbb{P}\left(\varphi_{s}^{-1}\left(X^{+} \right) > t \right) = 0.
\end{equation}
Now, it follows from (3.8), (3.10), Lemma 3.2 (with $Y = \varphi_{s}^{-1}\left(X^{+} \right)$, ~$p(t) = t$, and $h(t) = g(\log t)$, $t \geq 1$), and (3.9) that
\[
\begin{array}{ll}
& \mbox{$\displaystyle
\limsup_{n \to \infty} \frac{\log \left(n \mathbb{P}\left(X > s \sqrt{n g(\log n)} \right) \right)}{g(\log n)}$}\\
&\\
& \mbox{$\displaystyle
= \limsup_{n \to \infty} \frac{\log \left(n \mathbb{P}\left(\varphi_{s}^{-1}\left(X^{+} \right) > n \right) \right)}{g(\log n)}$}\\
& \\
& \mbox{$\displaystyle
= \limsup_{t \to \infty} \frac{\log \left(t \mathbb{P}\left(\varphi_{s}^{-1}\left(X^{+} \right) > t \right) \right)}{g(\log t)}$
~~(by (3.8), (3.10), and Lemma 3.2)}\\
& \\
& \mbox{$\displaystyle
= \limsup_{t \to \infty} \frac{\log \left(t \mathbb{P}\left(X > \varphi_{s}(t) \right) \right)}{g(\log t)}$}\\
& \\
& \mbox{$\displaystyle
= \limsup_{t \to \infty}
\frac{\log \left(\left(\varphi_{s}(t)\right)^{2} \mathbb{P}\left(X > \varphi_{s}(t) \right) \right)
- 2\log s - \log g(\log t)}{2^{\rho}g\left(\log \varphi_{s}(t) \right)}$
~~(by (3.9))}\\
& \\
& \mbox{$\displaystyle
= \limsup_{t \to \infty}
\frac{\log \left(\left(\varphi_{s}(t)\right)^{2}
\mathbb{P}\left(X > \varphi_{s}(t) \right) \right)}{2^{\rho}g\left(\log \varphi_{s}(t) \right)}$}\\
& \\
& \mbox{$\displaystyle
= \limsup_{x \to \infty}
\frac{\log \left(x^{2} \mathbb{P}(X > x) \right)}{2^{\rho}g(\log x)}$
~~(let $\displaystyle t = \varphi_{s}^{-1}(x)$)}\\
& \\
& \mbox{$\displaystyle
= - \frac{\overline{\lambda}_{1}}{2^{\rho}}$~~(by (2.1));}
\end{array}
\]
i.e., the first assertion of (3.6) holds. Similarly, the second assertion of (3.6) also follows from (3.8), (3.10), Lemma 3.2, and (3.9).

We now prove (3.7). Write $h(t) = g(\log t)$, $t \geq 1$. Under the given conditions, it is easy to see that
$h(\cdot):~[1, \infty) \rightarrow [0, \infty)$ is
a continuous and non-decreasing slowly varying function such that $\lim_{t \to \infty}h(t) = \infty$.
By the Karamata representation theorem (see, e.g., Bingham, Goldie, and Teugels [2, Theorem 1.3.1]), there exist two measurable functions
$c(\cdot)$ and $\varepsilon(\cdot)$ and a constant $d > 0$ such that
\begin{equation}
\left\{
\begin{array}{ll}
& \mbox{$\displaystyle
h(t) = c(t) \exp\left(\int_{d}^{t} \frac{\varepsilon(x)}{x}dx \right), ~t \geq d,$}\\
&\\
& \mbox{$\displaystyle
\lim_{t \to \infty} c(t) = c \in (0, \infty),
~\mbox{and}~ \lim_{t \to \infty} \varepsilon(t) = 0.$}
\end{array}
\right.
\end{equation}
Write
\[
\phi(t) = (1/c)\sqrt{t} \exp\left(-\int_{d}^{t} \frac{\varepsilon(x)}{x}dx \right), ~t \geq d.
\]
It follows from (3.11) that
\begin{equation}
\lim_{t \to \infty} \frac{\sqrt{t}/g(\log t)}{\phi(t)} =1,
\end{equation}
\begin{equation}
\lim_{t \to \infty} \frac{\phi(2t)}{\phi(t)} = \sqrt{2},
\end{equation}
\begin{equation}
\lim_{t \to \infty} \frac{\log \left(\frac{t}{\phi^{2}(t)}\right)}{g(\log t)} = \lim_{t \to \infty} \frac{2\log g(\log t)}{g(\log t)} = 0,
\end{equation}
\begin{equation}
\lim_{t \to \infty} \frac{g\left(\log \left(s\phi(t) \right)\right)}{g(\log t)}
= \lim_{t \to \infty} \frac{g\left(\frac{1}{2} \log t + \log s - \log g(\log t) \right)}{g(\log t)}
= \frac{1}{2^{\rho}}
\end{equation}
and
\begin{equation}
\phi^{\prime}(t) = \frac{\phi(t)}{t} \left(\frac{1}{2} - \varepsilon(t) \right)
> 0 ~\mbox{ultimately}.
\end{equation}
Clearly, from (3.12), we see that (3.7) is equivalent to, for all $s > 0$,
\begin{equation}
\left\{
\begin{array}{ll}
& \mbox{$\displaystyle
\limsup_{n \to \infty} \frac{\log \left(n \mathbb{P}(X > s \phi(n)) \right)}{g(\log n)}
= - \frac{\overline{\lambda}_{1}}{2^{\rho}},$}\\
&\\
& \mbox{$\displaystyle
\liminf_{n \to \infty} \frac{\log \left(n \mathbb{P} (X > s \phi(n)) \right)}{g(\log n)}
= - \frac{\underline{\lambda}_{1}}{2^{\rho}}.$}
\end{array}
\right.
\end{equation}
Note that (3.16) implies that there exists a constant $d_{1} > d$ such that $\phi(t)$ is
a continuous and strictly increasing function on $[d_{1}, \infty)$. Define
\[
\psi(t) =
\left \{
\begin{array}{ll}
\mbox{$\displaystyle \frac{\phi\left(d_{1} \right)}{d_{1}}t$} & \mbox{if $\displaystyle 0 \leq t < d_{1}$,}\\
&\\
\mbox{$\displaystyle \phi(t)$} & \mbox{if $\displaystyle  t \geq d_{1}$.}
\end{array}
\right.
\]
Then $\psi(t)$ is a continuous and strictly increasing function on $[0, \infty)$ $\psi(0) = 0$ and
$\lim_{t \to \infty} \psi(t) = \infty$. Let $\psi^{-1}(\cdot)$ be the inverse function of $\psi(\cdot)$.
Clearly, (3.13) ensures that there exists $t_{1} > d_{1}$ such that $0 < \phi(2t) \leq 2 \phi(t)$,
$t \geq t_{1}$; i.e.,
\begin{equation}
0 < \psi(2t) \leq 2 \psi(t), ~t \geq t_{1}.
\end{equation}
Since $\mathbb{E}X^{2} < \infty$, we have
\[
\limsup_{t \to \infty} \left(\psi(t)\right)^{2} \mathbb{P}\left(\psi^{-1}\left(\frac{X^{+}}{s}\right) > t \right)
= \limsup_{t \to \infty} \left(\psi(t)\right)^{2} \mathbb{P}\left( \frac{X}{s} > \psi(t) \right)
= 0.
\]
Thus
\begin{equation}
\lim_{t \to \infty} \left(\psi(t)\right)^{2} \mathbb{P}\left(\psi^{-1}\left(\frac{X^{+}}{s}\right) > t \right) = 0.
\end{equation}
Now, it follows from (3.14), the definition of $\psi(\cdot)$, (3.18), (3.19), (3.8), Lemma 3.2
(with $Y = \psi^{-1}(X^{+}/s)$, ~$p(t) = \psi^{2}(t), ~t > t_{1}$, and $h(t) = g(\log t)$, $t \geq 1$), (3.15),
and (2.1)) that
\[
\begin{array}{ll}
& \mbox{$\displaystyle
\limsup_{n \to \infty} \frac{\log \left(n \mathbb{P}\left(X > s \phi(n) \right) \right)}{g(\log n)}$}\\
&\\
& \mbox{$\displaystyle
= \limsup_{n \to \infty}
\frac{\log \left(\phi^{2}(n) \mathbb{P}\left(X > s \phi(n) \right)\right) + \log \left(\frac{n}{\phi^{2}(n)}\right)}{g(\log n)}$}\\
&\\
& \mbox{$\displaystyle
= \limsup_{n \to \infty}
\frac{\log \left(\phi^{2}(n) \mathbb{P}\left(X > s \phi(n) \right)\right)}{g(\log n)}$
~~(by (3.14))}\\
&\\
& \mbox{$\displaystyle
= \limsup_{n \to \infty}
\frac{\log \left(\psi^{2}(n) \mathbb{P}\left(X > s \psi(n) \right)\right)}{g(\log n)}$
~~(by the definition of $\psi(\cdot)$)}\\
&\\
& \mbox{$\displaystyle
= \limsup_{n \to \infty} \frac{\log \left(\psi^{2}(n) \mathbb{P}\left(\psi^{-1}\left(X^{+}/s \right) > n \right) \right)}{g(\log n)}$}\\
&\\
& \mbox{$\displaystyle
= \limsup_{t \to \infty} \frac{\log \left(\psi^{2}(t) \mathbb{P}\left(\psi^{-1}\left(X^{+}/s \right) > t \right) \right)}{g(\log t)}$
~~(by (3.18), (3.19), (3.8), and Lemma 3.2)}\\
&\\
& \mbox{$\displaystyle
= \limsup_{t \to \infty} \frac{\log \left(\psi^{2}(t) \mathbb{P}\left(X > s \psi(t) \right) \right)}{g(\log t)}$}\\
&\\
& \mbox{$\displaystyle
= \limsup_{t \to \infty}
\frac{\log \left(\left(s\psi(t)\right)^{2} \mathbb{P}\left(X > s \psi(t) \right) \right) - 2 \log s}{2^{\rho}g(\log \left(s\psi(t)\right))}$
~~(by (3.15) and the definition of $\psi(\cdot)$)}\\
& \\
& \mbox{$\displaystyle
= \limsup_{x \to \infty} \frac{\log \left(x^{2} \mathbb{P}(X > x) \right)}{2^{\rho}g(\log x)}$
~~(let $\displaystyle t = \psi^{-1}(x/s)$)}\\
& \\
& \mbox{$\displaystyle
= - \frac{\overline{\lambda}_{1}}{2^{\rho}}$~~(by (2.1));}
\end{array}
\]
i.e., the first assertion of (3.17) follows. Following the same argument, the second assertion of (3.17) follows.
The proof of Lemma 3.3 is complete. ~$\Box$

\vskip 0.2cm

\begin{lemma}
For each $n \geq 2$, let $\{X_{n,i}; 1 \leq i \leq n \}$ be i.i.d. real-valued random variables such that
\begin{equation}
\mathbb{E}X_{n, 1} = 0~~\mbox{and}~~\lim_{n \to \infty} \mathbb{E}X_{n,1}^{2} = \sigma^{2} \in (0, \infty).
\end{equation}
and, for some given $g(\cdot) \in \mathcal{V}_{\rho}$ for some $\rho \geq 0$, there exists a sequence of
positive constants $\left\{\tau_{n}; ~n \geq 2 \right\}$ such that
\begin{equation}
\lim_{n \to \infty} \tau_{n} = 0~~\mbox{and}
~~\left|X_{n,1} \right| \leq \tau_{n} \sqrt{\frac{n}{g(\log n)}}~~\mbox{almost surely (a.s.)}
\end{equation}
Then we have
\begin{equation}
\lim_{n \rightarrow \infty} \frac{\log \mathbb{P}\left(\sum_{i=1}^{n} X_{n,i} > r \sqrt{n g(\log n)} \right)}{g(\log n)}
= - \frac{r^{2}}{2\sigma^{2}}~~\mbox{for all}~ r > 0
\end{equation}
and similarly,
\[
\lim_{n \rightarrow \infty} \frac{\log \mathbb{P}\left(\sum_{i=1}^{n} X_{n,i} < - r \sqrt{n g(\log n)} \right)}{g(\log n)}
= - \frac{r^{2}}{2\sigma^{2}}~~\mbox{for all}~ r > 0
\]
and
\[
\lim_{n \rightarrow \infty} \frac{\log \mathbb{P}\left(\left|\sum_{i=1}^{n} X_{n,i} \right| > r \sqrt{n g(\log n)} \right)}{g(\log n)}
= - \frac{r^{2}}{2\sigma^{2}}~~\mbox{for all}~ r > 0.
\]
\end{lemma}

\vskip 0.2cm

\noindent {\it Proof}~~For $n \geq 2$ and fixed $r > 0$, write
\[
B_{n} = \mbox{Var}\left(\sum_{i=1}^{n} X_{n,i} \right), ~M_{n} = \tau_{n} \sqrt{\frac{n}{g(\log n)}},
~\mbox{and}~ ~x_{n}(r) = r \sqrt{n g(\log n)}, ~ n \geq 1.
\]
Since, for each $n \geq 2$, $\{X_{n,i}; 1 \leq i \leq n \}$ are i.i.d. real-valued random variables with (3.20), we have
\begin{equation}
\lim_{n \to \infty} \frac{B_{n}}{n \sigma^{2}} = \lim_{n \to \infty} \frac{n \mathbb{E}X_{n,1}^{2}}{n \sigma^{2}} = 1.
\end{equation}
Thus it follows from (3.21) and (3.22) that
\begin{equation}
\max_{1 \leq i \leq n} \left|X_{n,i} \right| \leq M_{n}~ \mbox{a.s.}, ~n \geq 2,
\end{equation}
\begin{equation}
0 < x_{n}(r) M_{n} = \tau_{n} rn \leq B_{n} ~\mbox{for all sufficiently large}~ n,
\end{equation}
\begin{equation}
\lim_{n \to \infty} \frac{x_{n}(r)M_{n}}{B_{n}} = \frac{r}{\sigma^{2}}\lim_{n \to \infty}
\left(\frac{n \sigma^{2}}{B_{n}}\right) \tau_{n} = 0,
\end{equation}
and
\begin{equation}
\frac{x_{n}^{2}(r)}{B_{n}} \sim \frac{r^{2}}{\sigma^{2}}
\left(\frac{n \sigma^{2}}{B_{n}} \right) g(\log n) \to \infty
~(\mbox{since} ~\lim_{t \to \infty} g(t) = \infty).
\end{equation}
By (3.24), (3.25), and Lemma 7.1 (which is one of the Kolmogorov exponential inequalities) of
Petrov [27, page 240], for all sufficiently large $n$ we have
\[
\begin{array}{lll}
\mbox{$\displaystyle
\mathbb{P}\left(\sum_{i=1}^{n} X_{n,i} > r \sqrt{n g(\log n)} \right)$}
& = &
\mbox{$\displaystyle
\mathbb{P}\left(\sum_{i=1}^{n} X_{n,i}  > x_{n}(r) \right)$}\\
&&\\
& \leq &
\mbox{$\displaystyle \exp \left\{- \frac{x_{n}^{2}(r)}{2B_{n}}\left(1 - \frac{x_{n}(r)M_{n}}{2B_{n}} \right) \right\}$}
\end{array}
\]
and hence, it follows from $\lim_{t \to \infty} g(t) = \infty$,
(3.26), (3.27), and (3.23) that
\begin{equation}
\begin{array}{lll}
\mbox{$\displaystyle
\limsup_{n \to \infty} \frac{\log \mathbb{P}\left(\sum_{i=1}^{n} X_{n,i} > r \sqrt{n g(\log n)} \right)}{g(\log n)}$}
& \leq &
\mbox{$\displaystyle
\limsup_{n \to \infty} \frac{- \frac{x_{n}^{2}(r)}{2B_{n}}\left(1 - \frac{x_{n}(r)M_{n}}{2B_{n}} \right)}{g(\log n)}$}\\
&&\\
& \leq &
\mbox{$\displaystyle
- \frac{r^{2}}{2 \sigma^{2}} \liminf_{n \to \infty} \frac{n \sigma^{2}}{B_{n}}$}\\
&&\\
& = &
\mbox{$\displaystyle - \frac{r^{2}}{2 \sigma^{2}}$.}
\end{array}
\end{equation}
Now by (3.26), (3.27), and Lemma 7.2 (which also is one of the Kolmogorov exponential inequalities) of
Petrov [27, page 241], for every fixed $0 < \epsilon < 1$ and all sufficiently large $n$ we have
\[
\begin{array}{lll}
\mbox{$\displaystyle
\mathbb{P}\left(\sum_{i=1}^{n} X_{n,i} > r \sqrt{n g(\log n)} \right)$}
& = &
\mbox{$\displaystyle
\mathbb{P}\left(\sum_{i=1}^{n} X_{n,i}  > x_{n}(r) \right)$}\\
&&\\
& \geq &
\mbox{$\displaystyle \exp \left\{- \frac{x_{n}^{2}(r)}{2B_{n}}(1 - \epsilon) \right\}$}
\end{array}
\]
and hence, it follows from (3.23) that
\[
\begin{array}{lll}
\mbox{$\displaystyle
\liminf_{n \to \infty} \frac{\log \mathbb{P}\left( \sum_{i=1}^{n} X_{n,i} > r \sqrt{n g(\log n)} \right)}{g(\log n)}$}
& \geq &
\mbox{$\displaystyle
\liminf_{n \to \infty} \frac{- \frac{x_{n}^{2}(r)}{2B_{n}}(1 - \epsilon)}{g(\log n)}$}\\
&&\\
& \geq &
\mbox{$\displaystyle
- \frac{(1 - \epsilon)r^{2}}{2 \sigma^{2}} \limsup_{n \to \infty} \frac{n \sigma^{2}}{B_{n}}$}\\
&&\\
& = &
\mbox{$\displaystyle - \frac{(1 - \epsilon)r^{2}}{2 \sigma^{2}}$.}
\end{array}
\]
Thus, letting $\epsilon \searrow 0$, we get
\begin{equation}
\liminf_{n \to \infty} \frac{\log \mathbb{P}\left(\sum_{i=1}^{n} X_{n,i} > r \sqrt{n g(\log n)} \right)}{g(\log n)}
\geq - \frac{r^{2}}{2 \sigma^{2}}.
\end{equation}
Clearly, (3.28) and (3.29) together ensure (3.22). This completes the proof of Lemma 3.4. ~$\Box$

\section{Proof of Theorem 2.1}

We first state two basic facts which will be used in the proof of Theorem 2.1. Let $\left\{a_{n};~n \geq 1 \right\}$
and $\left\{b_{n};~n \geq 1 \right\}$ be two sequences
of real numbers. Then
\[
\limsup_{n \to \infty} \left(a_{n}\vee b_{n} \right)
= \left(\limsup_{n \to \infty} a_{n} \right) \vee \left(\limsup_{n \to \infty} b_{n} \right)
~~\mbox{and}~~\liminf_{n \to \infty} \left(a_{n}\vee b_{n} \right)
\leq \left(\limsup_{n \to \infty} a_{n} \right) \vee \left(\liminf_{n \to \infty} b_{n} \right).
\]

\vskip 0.2cm

\noindent {\it Proof of Theorem 2.1}~~Since $g(\cdot) \in \mathcal{V}_{\rho}$, we have
\[
\lim_{t \to \infty} \frac{g(\log(t \pm \mu ))}{g(\log t)} = 1
\]
and hence, in view of (2.1), (2.2), and (2.3), without of loss generality, we can assume that $\mu = 0$.

\vskip 0.2cm

{\bf (i)}~~We first give the proof of (2.4).

{\bf The lower bound part of (2.4)}.~~We first show that, for all $x > 0$,
\begin{equation}
\left \{
\begin{array}{ll}
& \mbox{$\displaystyle
- \frac{\overline{\lambda}_{1}}{2^{\rho}}
\leq \limsup_{n \to \infty} \frac{\log \mathbb{P}\left(S_{n} > x \sqrt{n g(\log n)} \right)}{g(\log n)},$}\\
& \\
& \mbox{$\displaystyle
- \frac{\underline{\lambda}_{1}}{2^{\rho}}
\leq \liminf_{n \to \infty} \frac{\log \mathbb{P}\left(S_{n} > x \sqrt{n g(\log n)} \right)}{g(\log n)}.$}
\end{array}
\right.
\end{equation}
Recall that $m(Y)$ is a median for a real-valued random variable $Y$.
Write $S_{0} = 0$ and
\[
m_{n} = \min_{0 \leq k \leq n} m\left(S_{k} \right), ~n \geq 1,
\]
Since $X_{1}, ..., X_{n}$ are i.i.d. random variables, we have
\[
\min_{0 \leq k \leq n} m\left(S_{n} - S_{k} \right) = m_{n},~~n \geq 1.
\]
Under the given conditions of Theorem 2.1, we have
\[
\frac{S_{n}}{\sqrt{n g(\log n)}} \to_{\mathbb{P}} 0,
\]
where ``$\to_{\mathbb{P}}$" stands for convergence in probability. Hence,
\begin{equation}
\lim_{n \rightarrow \infty} \frac{m_{n}}{\sqrt{n g(\log n)}} = \lim_{n \rightarrow \infty}
\frac{\min_{0 \leq k \leq n} m\left(S_{k} \right)}{\sqrt{n g(\log n)}} = 0.
\end{equation}
Thus, for any given $x > 0$, by (4.2) and Lemma 3.1, we have for all sufficiently large $n$,
\begin{equation}
\begin{array}{lll}
\mbox{$\displaystyle
    \mathbb{P} \left(\max_{1 \leq k \leq n} X_{k} > 2x \sqrt{n g(\log n)} \right)$}
& \leq &
\mbox{$\displaystyle
    \mathbb{P} \left(\max_{1 \leq k \leq n} X_{k} > x \sqrt{n g(\log n)} - 2m_{n} \right)$ ~(by (4.2))}\\
&&\\
& = &
\mbox{$\displaystyle
\mathbb{P} \left(\max_{1 \leq k \leq n} X_{k} +
\min_{0 \leq k \leq n} m\left(S_{k} \right) > x \sqrt{n g(\log n)} - m_{n} \right)$}\\
&&\\
& \leq &
\mbox{$\displaystyle
    \mathbb{P} \left(\max_{1 \leq k \leq n} \left(X_{k}
    + m\left(S_{k-1}\right) \right) > x \sqrt{n g(\log n)} - m_{n} \right)$}\\
&&\\
& \leq &
\mbox{$\displaystyle
    2 \mathbb{P} \left(\max_{1 \leq k \leq n} S_{k} > x \sqrt{n g(\log n)} - m_{n} \right)$~~(by Lemma 3.1)}\\
&&\\
& \leq &
\mbox{$\displaystyle
    2 \mathbb{P} \left(\max_{1 \leq k \leq n} \left(S_{k}
    + m\left(S_{n} - S_{k}\right) \right) > x \sqrt{n g(\log n)} \right)$}\\
&&\\
& \leq &
\mbox{$\displaystyle
    4 \mathbb{P}\left(S_{n} > x \sqrt{n g(\log n)} \right)$~~(by Lemma 3.1).}
\end{array}
\end{equation}
Using the method used in the proof of Lemma 3.4 of Li and Miao [19], we get
\begin{equation}
\frac{1 \wedge \left(n \mathbb{P}\left(X > 2x \sqrt{ng(\log n)} \right) \right)}{2}
\leq \mathbb{P}\left(\max_{1 \leq k \leq n} X_{k} > 2x\sqrt{ng(\log n)} \right),
~n \geq 1.
\end{equation}
Note that $\mathbb{E}X^{2} < \infty$ and $\lim_{t \to \infty} g(t) = \infty$ ensure that
\[
\lim_{n \to \infty} n \mathbb{P}\left(X > 2x \sqrt{ng(\log n)} \right) = 0.
\]
Thus, under the given conditions of $g(\cdot)$, it follows from (4.3), (4.4), and Lemma 3.3 that
\[
\begin{array}{lll}
\mbox{$\displaystyle
\limsup_{n \to \infty} \frac{\log \mathbb{P}\left(S_{n} > x \sqrt{n g(\log n)} \right)}{g(\log n)}$}
& \geq &
\mbox{$\displaystyle
\limsup_{n \to \infty} \frac{\log \left(\frac{n \mathbb{P}\left(X > 2x \sqrt{n g(\log n)} \right)}{8} \right)}{g(\log n)}$}\\
&&\\
& = &
\mbox{$\displaystyle
\limsup_{n \to \infty} \frac{\log \left(n \mathbb{P}\left(X > 2x \sqrt{n g(\log n)} \right) \right)}{g(\log n)}$}\\
&&\\
& = &
\mbox{$\displaystyle
- \frac{\overline{\lambda}_{1}}{2^{\rho}} ~~\mbox{for all}~ x > 0$;}
\end{array}
\]
i.e., the first half of (4.1) holds. In the same vein, the second half of (4.1) follows.

In the following, we aim to establish the inequality
\begin{equation}
- \frac{x^{2}}{2 \sigma^{2}}
\leq \liminf_{n \to \infty} \frac{\log \mathbb{P}\left(S_{n} > x \sqrt{n g(\log n)} \right)}{g(\log n)}
~~\mbox{for all}~ x > 0.
\end{equation}
Since $\mathbb{E}X^{2} < \infty$, we have
\[
\lim_{t \to \infty}
\frac{1}{\delta^2}\mathbb{E}\left(X^{2}I\{|X| > \delta t \} \right) = 0~~\mbox{for all}~ \delta > 0.
\]
Thus there exists a sequence of positive constants
$\left \{\delta_{n};~ n \geq 1 \right \}$ with $\delta_{n} \searrow 0$ as $n \to \infty$
such that
\[
\lim_{n \to
\infty} \frac{1}{\delta_n^2}\mathbb{E}\left(|X|^2I\left\{|X| >
\delta_{n} \sqrt{\frac{n}{g(\log n)}} \right \} \right) = 0,
\]
which implies
\begin{equation}
\lim_{n \to
\infty} \sqrt{\frac{n}{g(\log n)}} \mathbb{E}\left(|X|I\left\{|X| >
\delta_{n} \sqrt{\frac{n}{g(\log n)}} \right \} \right) = 0
\end{equation}
and
\begin{equation}
\lim_{n \to
\infty}\frac{n}{g(\log n)}\mathbb{P}\left(|X| >
\delta_{n} \sqrt{\frac{n}{g(\log n)}}\right)=0
\end{equation}
in view of the following estimates, for $t > 0$,
\[
t\mathbb{E}\left(|X|I\{|X| > \delta t \}\right) \leq \frac{1}{\delta}
\mathbb{E}\left(X^{2}I\{|X| > \delta t \} \right)
~~\mbox{and}~~
t^{2}\mathbb{P}\left(|X|>\delta t\right)\le \frac{1}{\delta^2}\mathbb{E}\left(|X|^2I\left\{|X| >
\delta t\right\}\right).
\]
Write, for $n \geq 2$,
\begin{equation}
\hat{\delta}_{n} = \delta_{n} \vee \frac{1}{\sqrt{g(\log n)}}, ~c_{n} = \hat{\delta}_{n} \sqrt{\frac{n}{g(\log n)}},
~\mu_{n} = \mathbb{E}\left(XI\left\{|X| \leq c_{n} \right \} \right), ~\mbox{and}~ p_{n} = \mathbb{P}\left(|X| > c_{n} \right).
\end{equation}
It follows from (4.6) (since $c_{n} \geq \delta_{n} \sqrt{\frac{n}{g(\log n)}}$,~$n \geq 2$)  and (4.7) that
\begin{equation}
\lim_{n \to \infty}\sqrt{\frac{n}{g(\log n)}} \mathbb{E}\left(|X|I\left\{|X| >
c_{n}\right\} \right) = 0~\mbox{ and }~\lim_{n \to \infty}\frac{np_{n}}{g(\log n)}= 0.
\end{equation}
Since $\mu = \mathbb{E}X = 0$, we have
\[
\frac{n\left|\mu_{n} \right|}{\sqrt{ng(\log n)}} \leq \sqrt{\frac{n}{g(\log n)}} \mathbb{E}\left(|X|I\left\{|X| >
c_{n} \right \} \right), ~n \geq 2
\]
and hence, from the first half of (4.9),
\begin{equation}
\lim_{n\to\infty}\frac{n\mu_n}{\sqrt{ng(\log n)}}=0.
\end{equation}
Since $c_{n} \geq \frac{\sqrt{n}}{g(\log n)} \to \infty$ as $n \to \infty$, we conclude that $p_{n} \to 0$ as $n \to \infty$ and
hence, from the second half of (4.9),
\begin{equation}
\lim_{n \to \infty}\frac{\log(1-p_{n})^{n}}{g(\log n)}
= \lim_{n\to\infty}\frac{n\log(1-p_{n})}{g(\log n)}
= \lim_{n\to\infty}\frac{-np_{n}}{g(\log n)} = 0.
\end{equation}
Denote the conditional distribution of $X$ given $|X|\leq c_{n}$ as $F_{n}$; that is,
\[
F_{n}(x)=\mathbb{P}\left(X\leq x \big||X|\leq c_{n} \right) = \left\{
 \begin{array}{ll}
 0 & \mbox{if $\displaystyle x<-c_{n},$} \\
 &\\
 \mbox{$\displaystyle \frac{1}{1-p_{n}}\mathbb{P}\Big(-c_{n} \le X\leq x\Big)$}
 & \mbox{if $\displaystyle -c_{n} \leq x \leq c_{n},$}  \\
 &\\
 1 & \mbox{if $\displaystyle x > c_{n}.$}
 \end{array}
 \right.
\]
Conditional on $\{\max_{1\le k \leq n} \left|X_{k} \right|\le c_n\}$,  $X_{1}, \cdots, X_{n}$ are i.i.d. random variables
with distribution function $F_{n}$. Let $\widetilde{X}_{n,1},\cdots, \widetilde{X}_{n,n}$ be independent random variables
with common distribution $F_{n}$. Since $\lim_{n \to \infty} p_{n} = 0$, it is easy to see that
\begin{equation}
\tilde{\mu}_{n}:=\mathbb{E} \Big(\widetilde{X}_{n,1}\Big) = \frac{\mu_{n}}{1-p_{n}} \to 0
~~\mbox{as}~ n \to \infty,
\end{equation}
and
\[
\tilde{\sigma}_{n}^{2} = \mathbb{E}\Big(\widetilde{X}_{n,1} - \tilde{\mu}_{n} \Big)^{2}
= \frac{1}{1-p_{n}}\mathbb{E}\left(X^{2}I\left\{|X| \leq c_{n} \right \}\right) - \tilde{\mu}_{n}^{2}
\to \sigma^{2}.
\]
Clearly, $\big\{U_{n,i}:=\widetilde{X}_{n,i}-\tilde{\mu}_{n}; ~1\leq i\leq n \big\}$ are i.i.d. random variables with mean $0$
and variance $\tilde{\sigma}_{n}^{2} \to \sigma^{2}$ as $n \to \infty$ and bounded by $2c_{n}= \tau_{n} \sqrt{\frac{n}{g(\log n)}}$,
$n \geq 2$ where $\tau_{n} = 2 \hat{\delta}_{n} \to 0$ as $n \to \infty$. Thus, replacing $\big\{X_{n,i}; ~1 \leq i \leq n, n \geq 2 \big\}$
with $\big\{U_{n,i}; ~1 \leq i \leq n, n \geq 2 \big\}$, all conditions of Lemma 3.4 are satisfied and hence, (3.22) holds for
$\big\{U_{n,i}; ~1 \leq i \leq n, n \geq 2 \big\}$. Thus, for any fixed $x>0$ and $\epsilon > 0$, we have
\[
\begin{array}{ll}
& \mbox{$\displaystyle \mathbb{P}\left(S_{n} > x\sqrt{ng(\log n)}\right)$}\\
&\\
& \mbox{$\displaystyle
\geq \mathbb{P}\left(\big\{S_{n} > x\sqrt{ng(\log n)} \big\} \bigcap
\big\{\max_{1\leq k \leq n}|X_{k}|\leq c_{n} \big\} \right)$}\\
&\\
&
\mbox{$\displaystyle
 = \mathbb{P}\left(\big\{\sum_{i=1}^{n} X_{i}I\left\{|X_{i}| \leq c_{n} \right\} > x\sqrt{ng(\log n)} \big\} \bigcap
\big\{\max_{1\leq k \leq n}|X_{k}|\leq c_{n} \big\} \right)$}\\
&\\
& \mbox{$\displaystyle
= \mathbb{P}\left(\sum_{i=1}^{n} X_{i}I\left\{|X_{i}| \leq c_{n} \right\} > x\sqrt{ng(\log n)}
\Big|\max_{1\leq k \leq n}|X_{k}|\leq c_{n}\right)\mathbb{P}
\left(\max_{1\leq k \leq n}|X_{k}|\le c_{n}\right)$}\\
&\\
& \mbox{$\displaystyle
= \mathbb{P}\left(\sum_{i=1}^{n}\widetilde{X}_{n,i} > x\sqrt{ng(\log n)} \right)\left(1-p_{n} \right)^{n}$}\\
&\\
& \mbox{$\displaystyle
\geq \mathbb{P}\left(\sum_{i=1}^{n}\left(\widetilde{X}_{n,i}-\tilde{\mu}_{n} \right) > x\sqrt{ng(\log n)}
+n\left|\tilde{\mu}_{n} \right| \right)\left(1-p_{n} \right)^{n}$}\\
&\\
& \mbox{$\displaystyle
\geq \mathbb{P}\left(\sum^n_{i=1}U_{n,i} > (x+\epsilon)\sqrt{ng(\log n)}\right)\left(1-p_{n} \right)^{n}$}
\end{array}
\]
for all sufficiently large $n$. In the last step, we have used (4.10) and (4.12) to conclude that
$n|\tilde{\mu}_n|\leq \epsilon\sqrt{ng(\log n)}$ for all sufficiently large $n$. Therefore, from
Lemma 3.4 and (4.11) we obtain that
\[
\begin{array}{ll}
& \mbox{$\displaystyle
\liminf_{n \to \infty} \frac{\log \mathbb{P}\left(S_{n} > x\sqrt{ng(\log n)}\right)}{g(\log n)}$}\\
& \\
& \mbox{$\displaystyle
\geq \lim_{n \to \infty}\frac{\log\mathbb{P}\left(\sum_{i=1}^{n} U_{n,i}
> (x + \epsilon)\sqrt{ng(\log n)}\right)}{g(\log n)}
+ \lim_{n\to\infty}\frac{\log\Big(\left(1-p_{n} \right)^{n}\Big)}{g(\log n)}$}\\
& \\
& \mbox{$\displaystyle
= -\frac{(x + \epsilon)^2}{2\sigma^{2}},$}
\end{array}
\]
which yields (4.5) by letting $\epsilon \searrow 0$.

Clearly, (4.1) and (4.5) together ensure the following lower bound of (2.4):
\begin{equation}
\left\{
\begin{array}{ll}
&
\mbox{$\displaystyle - \left(\frac{x^{2}}{2\sigma^{2}} \wedge \frac{\overline{\lambda}_{1}}{2^{\rho}} \right)
\leq \limsup_{n \rightarrow \infty} \frac{\log \mathbb{P}\left(S_{n} > x \sqrt{n g(\log n)} \right)}{g(\log n)}
~\mbox{for all}~x > 0$,}\\
&\\
&
\mbox{$\displaystyle - \left(\frac{x^{2}}{2\sigma^{2}} \wedge \frac{\underline{\lambda}_{1}}{2^{\rho}} \right)
\leq \liminf_{n \rightarrow \infty} \frac{\log \mathbb{P}\left(S_{n} > x \sqrt{n g(\log n)} \right)}{g(\log n)}
~\mbox{for all}~x > 0$.}
\end{array}
\right.
\end{equation}

\vskip 0.2cm

{\bf The upper bound part of (2.4)}.~~To complete the proof of (2.4) we have to establish the
following upper bound of (2.4):
\begin{equation}
\left\{
\begin{array}{ll}
&
\mbox{$\displaystyle
\limsup_{n \rightarrow \infty} \frac{\log \mathbb{P}\left(S_{n} > x \sqrt{n g(\log n)} \right)}{g(\log n)}
\leq - \left(\frac{x^{2}}{2\sigma^{2}} \wedge \frac{\overline{\lambda}_{1}}{2^{\rho}} \right)
~\mbox{for all}~x > 0$,}\\
&\\
&
\mbox{$\displaystyle
\liminf_{n \rightarrow \infty} \frac{\log \mathbb{P}\left(S_{n} > x \sqrt{n g(\log n)} \right)}{g(\log n)}
\leq - \left(\frac{x^{2}}{2\sigma^{2}} \wedge \frac{\underline{\lambda}_{1}}{2^{\rho}} \right)
~\mbox{for all}~x > 0$.}
\end{array}
\right.
\end{equation}
Write, for $n \geq 2$
\[
V_{n,i} = X_{i}I\left\{\left|X_{i}\right| \leq c_{n} \right \} - \mu_{n}, ~ Y_{n,i} =
X_{i}I\left\{X_{i} > c_{n} \right\} + \mu_{n}, ~i = 1, 2, ..., n,
\]
where $\big\{c_{n}; n\geq 2 \big\}$ and $\big\{\mu_{n}; n \geq 2 \big\}$ are defined in (4.8).
Clearly, $\big\{V_{n,i}; ~1\leq i\leq n \big\}$ are i.i.d. random variables with mean $0$
and variance $\mbox{Var}\big(V_{n,1} \big) \to \sigma^{2}$ as $n \to \infty$ and bounded by $2c_{n}= \tau_{n} \sqrt{\frac{n}{g(log n)}}$,
$n \geq 2$ where $\tau_{n} = 2 \hat{\delta}_{n} \to 0$ as $n \to \infty$. Thus, replacing $\big\{X_{n,i}; ~1 \leq i \leq n, n \geq 2 \big\}$
with $\big\{V_{n,i}; ~1 \leq i \leq n, n \geq 2 \big\}$, all conditions of Lemma 3.4 are satisfied and hence, (3.22) holds for
$\big\{V_{n,i}; ~1 \leq i \leq n, n \geq 2 \big\}$. Note that, for $n \geq 2$,
\[
S_{n} \leq \sum_{i=1}^{n} X_{i} I\left\{X > -c_{n} \right\} = \sum_{i=1}^{n} V_{n,i} + \sum_{i=1}^{n} Y_{n,i}
~~\mbox{and}~~
\frac{\sqrt{n}}{g(\log n)} \leq c_{n}.
\]
Thus, for any fixed $x>0$ and $0 < \epsilon < x$, it follows from (4.10) that, for all sufficiently large $n$,
\begin{equation}
\begin{array}{ll}
& \mbox{$\displaystyle
\left \{ S_{n} > x \sqrt{n g(\log n)} \right\}$}\\
&\\
& \mbox{$\displaystyle
\subseteq \left \{\sum_{i=1}^{n} V_{n,i} + \sum_{i=1}^{n} Y_{n,i} > x \sqrt{n g(\log n)} \right\}$}\\
&\\
& \mbox{$\displaystyle
\subseteq \left \{ \sum_{i=1}^{n} V_{n,i}  > (x - \epsilon) \sqrt{n g(\log n)} \right\}
\bigcup \left \{ \sum_{i=1}^{n} X_{i}I\left\{X_{i} > c_{n} \right \}
 + n \mu_{n} > \epsilon \sqrt{n g(\log n)} \right\}$}\\
&\\
& \mbox{$\displaystyle
\subseteq \left \{ \sum_{i=1}^{n} V_{n,i} > (x - \epsilon) \sqrt{n g(\log n)} \right\}
\bigcup \left \{ \max_{1 \leq i \leq n} X_{i} > \frac{\sqrt{n}}{g(\log n)} \right\}$.}
\end{array}
\end{equation}
Then it follows from (4.15) and Lemmas 3.3 and 3.4 that
\[
\begin{array}{ll}
& \mbox{$\displaystyle
\limsup_{n \rightarrow \infty}
\frac{\log \mathbb{P}\left( S_{n} > x \sqrt{n g(\log n)} \right)}{g(\log n)}$}\\
& \\
& \mbox{$\displaystyle
\leq \limsup_{n \rightarrow \infty}
\frac{\log \left(\mathbb{P}\left( \sum_{i=1}^{n} V_{n,i} > (x - \epsilon) \sqrt{n g(\log n)} \right)
+ \mathbb{P}\left(\max_{1 \leq i \leq n} X_{i} > \frac{\sqrt{n}}{g(\log n)} \right) \right)}{g(\log n)}$
~~(by (4.15))}\\
& \\
& \mbox{$\displaystyle
\leq \left(\limsup_{n \rightarrow \infty}
\frac{\log \mathbb{P}\left(\sum_{i=1}^{n} V_{n,i} > (x - \epsilon) \sqrt{n g(\log n)}\right)}{g(\log n)}\right)
\vee \left(\limsup_{n \rightarrow \infty}
\frac{\log \left(n \mathbb{P}\left(X > \frac{\sqrt{n}}{g(\log n)}\right) \right)}{g(\log n)} \right)$}\\
&\\
& \mbox{$\displaystyle = \left( - \frac{(x - \epsilon)^{2}}{2
\sigma^{2}} \right) \vee \left(- \frac{\overline{\lambda}_{1}}{2^{\rho}}
\right)$ ~~(by (3.22) and the first half of (3.7)).}
\end{array}
\]
and
\[
\begin{array}{ll}
& \mbox{$\displaystyle
\liminf_{n \rightarrow \infty}
\frac{\log \mathbb{P}\left( S_{n} > x \sqrt{n g(\log n)} \right)}{g(\log n)}$}\\
& \\
& \mbox{$\displaystyle
\leq \liminf_{n \rightarrow \infty}
\frac{\log \left(\mathbb{P}\left( \sum_{i=1}^{n} V_{n,i} > (x - \epsilon) \sqrt{n g(\log n)} \right)
+ \mathbb{P}\left(\max_{1 \leq i \leq n} X_{i} > \frac{\sqrt{n}}{g(\log n)} \right) \right)}{g(\log n)}$
~~(by (4.15))}\\
& \\
& \mbox{$\displaystyle
\leq \left(\limsup_{n \rightarrow \infty}
\frac{\log \mathbb{P}\left(\sum_{i=1}^{n} V_{n,i} > (x - \epsilon) \sqrt{n g(\log n)}\right)}{g(\log n)}\right)
\vee \left(\liminf_{n \rightarrow \infty}
\frac{\log \left(n \mathbb{P}\left(X > \frac{\sqrt{n}}{g(\log n)}\right) \right)}{g(\log n)} \right)$}\\
&\\
& \mbox{$\displaystyle = \left( - \frac{(x - \epsilon)^{2}}{2
\sigma^{2}} \right) \vee \left(- \frac{\underline{\lambda}_{1}}{2^{\rho}}
\right)$ ~~(by (3.22) and the second half of (3.7)}
\end{array}
\]
which yields (4.14) by letting $\epsilon \searrow 0$.

\vskip 0.2cm

{\bf (ii)}~~Clearly, replacing $\left\{X, X_{n};~n \geq 1 \right\}$ with $\left\{-X, -X_{n};~n \geq 1 \right\}$,
(2.5) follows from (2.4) and (2.2).

\vskip 0.2cm

{\bf (iii)}~~To complete the proof of Theorem 2.1 we now establish (2.6). Since $\sigma^{2} < \infty$ and $\mu = 0$,
replacing $\left\{X, X_{n};~n \geq 1 \right\}$  with $\left\{-X, -X_{n};~n \geq 1 \right\}$, it follows
from (4.3) that, for any given $x > 0$,
\[
\mathbb{P}\left(\max_{1 \leq k \leq n}(-X_{k}) > 2x \sqrt{n g(\log n)} \right)
\leq 4 \mathbb{P} \left(-S_{n} > x \sqrt{n g(\log n)} \right)
~~\mbox{for all sufficiently large}~n.
\]
Together with (4.3) this implies that
\begin{equation}
\mathbb{P}\left(\max_{1 \leq k \leq n}\left|X_{k}\right| > 2x \sqrt{n g(\log n)} \right)
\leq 4 \mathbb{P} \left(\left|S_{n}\right| > x \sqrt{n g(\log n)} \right)
~~\mbox{for all sufficiently large}~n.
\end{equation}
In the same vein as establishing (4.1), it follows from (4.16) and (3.6) (i.e., the first part
of Lemma 3.3) that
\begin{equation}
\left \{
\begin{array}{ll}
& \mbox{$\displaystyle
- \frac{\overline{\lambda}}{2^{\rho}}
\leq \limsup_{n \to \infty} \frac{\log \mathbb{P}\left(\left|S_{n}\right| > x \sqrt{n g(\log n)} \right)}{g(\log n)},$}\\
& \\
& \mbox{$\displaystyle
- \frac{\underline{\lambda}}{2^{\rho}}
\leq \liminf_{n \to \infty} \frac{\log \mathbb{P}\left(\left|S_{n}\right| > x \sqrt{n g(\log n)} \right)}{g(\log n)}.$}
\end{array}
\right.
\end{equation}
Since, for $n \geq 1$ and $x > 0$,
$\mathbb{P}\left(S_{n} > x \sqrt{n g(\log n)} \right) \leq \mathbb{P}\left(\left|S_{n}\right| > x \sqrt{n g(\log n)} \right)$,
(4.5) and (4.17) together ensure the following lower bound of (2.6):
\begin{equation}
\left\{
\begin{array}{ll}
&
\mbox{$\displaystyle - \left(\frac{x^{2}}{2\sigma^{2}} \wedge \frac{\overline{\lambda}}{2^{\rho}} \right)
\leq \limsup_{n \rightarrow \infty} \frac{\log \mathbb{P}\left(\left|S_{n} \right| > x \sqrt{n g(\log n)} \right)}{g(\log n)}
~\mbox{for all}~x > 0$,}\\
&\\
&
\mbox{$\displaystyle - \left(\frac{x^{2}}{2\sigma^{2}} \wedge \frac{\underline{\lambda}}{2^{\rho}} \right)
\leq \liminf_{n \rightarrow \infty} \frac{\log \mathbb{P}\left(\left|S_{n}\right| > x \sqrt{n g(\log n)} \right)}{g(\log n)}
~\mbox{for all}~x > 0$.}
\end{array}
\right.
\end{equation}
Now note that, for $n \geq 2$,
\[
\begin{array}{lll}
\mbox{$\displaystyle \left|S_{n} \right|$}
& = &
\mbox{$\displaystyle
\left|\sum_{i=1}^{n} \left(X_{i}I\left\{|X_{i}| \leq c_{n} \right\} - \mu_{n} \right) +
\sum_{i=1}^{n} X_{i}I\left\{|X_{i}| > c_{n} \right\} + n \mu_{n} \right|$}\\
&&\\
& \leq &
\mbox{$\displaystyle \left|\sum_{i=1}^{n}V_{n,i} \right|
+ \sum_{i=1}^{n} \left|X_{i} \right|I\left\{|X_{i}| > c_{n} \right\} + \left|n \mu_{n} \right|$},
\end{array}
\]
where $\big\{c_{n}; n\geq 2 \big\}$, $\big\{\mu_{n}; n \geq 2 \big\}$, and $\big\{V_{n,i}; ~1 \leq i \leq n, n \geq 2 \big\}$
are the same as in (4.15). Thus, for any fixed $x>0$ and $0 < \epsilon < x$, by the same argument as in (4.15),
it follows from (4.10) that, for all sufficiently large $n$,
\begin{equation}
\left \{ \left|S_{n} \right| > x \sqrt{n g(\log n)} \right\}
\subseteq \left \{ \left|\sum_{i=1}^{n} V_{n,i} \right| > (x - \epsilon) \sqrt{n g(\log n)} \right\}
\bigcup \left \{ \max_{1 \leq i \leq n} \left|X_{i} \right| > \frac{\sqrt{n}}{g(\log n)} \right\}
\end{equation}
Then, for any fixed $x > 0$ and $0 < \epsilon < x$, it follows from (4.19), the second part of Lemma 3.3 (i.e., (3.7)
with $X$, $\overline{\lambda}_{1}$, and $\underline{\lambda}_{1}$ replaced by $|X|$, $\overline{\lambda}$, and
$\underline{\lambda}$ respectively) and Lemma 3.4 (with $\big\{X_{n,i}; ~1 \leq i \leq n, n \geq 2 \big\}$
replaced by $\big\{V_{n,i}; ~1 \leq i \leq n, n \geq 2 \big\}$) that
\begin{equation}
\left\{
\begin{array}{ll}
&
\mbox{$\displaystyle
\limsup_{n \rightarrow \infty} \frac{\log \mathbb{P}\left(\left|S_{n} \right| > x \sqrt{n g(\log n)} \right)}{g(\log n)}
\leq  - \left(\frac{(x-\epsilon)^{2}}{2\sigma^{2}} \wedge \frac{\overline{\lambda}}{2^{\rho}} \right)$,}\\
&\\
&
\mbox{$\displaystyle
\liminf_{n \rightarrow \infty} \frac{\log \mathbb{P}\left(\left|S_{n}\right| > x \sqrt{n g(\log n)} \right)}{g(\log n)}
\leq - \left(\frac{(x - \epsilon)^{2}}{2\sigma^{2}} \wedge \frac{\underline{\lambda}}{2^{\rho}} \right)$.}
\end{array}
\right.
\end{equation}
Clearly, (2.6) follows from (4.18) and (4.20) by letting $\epsilon \searrow 0$. This completes the proof of Theorem 2.1.
~$\Box$

\section{Proof of Theorem 2.2}

In this section, we give the proof of Theorem 2.2.

\vskip 0.2cm

\noindent {\it Proof of Theorem 2.2}~~~~Since $g(\cdot) \in \mathcal{V}_{\rho}$, we have
\[
\lim_{t \to \infty} \frac{g(\log(t \pm \eta ))}{g(\log t)} = 1
\]
and hence, in view of (2.1), (2.2), and (2.3), without of loss generality, we can assume that $\eta = 0$.

We only give the proof of Theorem 2.2 (i), the proofs of Theorem 2.2 (ii) and (iii)
are left to the reader.

We first establish the implication (2.9) $\Rightarrow$ (2.7). Since $\mathbb{E}|X| = \infty$ ensures that $\mathbb{E}X^{2} = \infty$,
(2.9) implies the following three cases to consider:

\noindent {\it Case I}~~~$\mathbb{E}|X| < \infty$ and $\mathbb{E}X \neq 0$;

\noindent {\it Case II}~~$\mathbb{E}X^{2} = \infty$;

\noindent {\it Case III}~$\mathbb{E}X = 0$, $\mathbb{E}X^{2} = \sigma^{2} \in (0, \infty)$, and $\overline{\lambda} = \underline{\lambda} = 0$.

For Case I, by the law of large numbers, for all $x$,
\[
\lim_{n \to \infty} \mathbb{P}\left(\left|S_{n} \right| > x \sqrt{n g(\log n)} \right)
= \lim_{n \to \infty} \mathbb{P}\left(\left|\frac{S_{n}}{n} \right| > x \sqrt{\frac{g(\log n)}{n}} \right) = 1
\]
which ensures (2.7).

For Case II, we consider $\tilde{S}_{n} = \sum_{i=1}^{n} Y_{i}$, $n \geq 1$, where $\left\{Y, Y_{n}; n \geq 1 \right \}
= \left \{X - X^{\prime}, X_{n} - X^{\prime}_{n}; n \geq 1 \right \}$, $\{X^{\prime},~X_{n}^{\prime };~n\geq 1\}$ is
an independent copy of $\{X,~X_{n};~n\geq 1\}$. Note that $\left \{Y, Y_{n}; ~n \geq 1 \right \}$ is
a sequence of i.i.d. symmetric real-valued random variables with $\mathbb{E}Y^{2} = \infty$ (since $\mathbb{E}X^{2} = \infty$).
Thus, by the second part of Lemma 6.5 of Ledoux and Talagrand [18, pages 153-154], we have, for each $n \geq 1$,
\begin{equation}
\begin{array}{lll}
\mbox{$\displaystyle
\mathbb{P}\left(\left|\sum_{i=1}^{n} Y_{i}(c) \right| > 2x\sqrt{ng(\log n)} \right)$}
& \leq &
\mbox{$\displaystyle
2 \mathbb{P}\left(\left|\tilde{S}_{n} \right| > 2x\sqrt{ng(\log n)} \right)$}\\
&&\\
&\leq& \mbox{$\displaystyle
4 \mathbb{P}\left(\left|S_{n} \right| > x \sqrt{ng(\log n)} \right)~~\mbox{for all}~x > 0,$}\\
\end{array}
\end{equation}
where $Y(c) = YI\{|Y| \leq c \}$, and $Y_{n}(c) = Y_{n} I\left\{|Y_{n}| \leq c \right \}$, $n \geq 1$ and $c > 0$
is large enough such that $\mathbb{E}Y^{2}(c) > 0$. Then, by Lemma 3.4 and (5.1), we have
\begin{equation}
\begin{array}{lll}
\mbox{$\displaystyle
- \frac{(2x)^{2}}{2 \mathbb{E}Y^{2}(c)}$}
&=&
\mbox{$\displaystyle \liminf_{n \to \infty}
\frac{\log \mathbb{P}\left(\left|\sum_{i=1}^{n} Y_{i}(c) \right| > 2x\sqrt{ng(\log n)} \right)}{g(\log n)}$ ~(by Lemma 3.4)}\\
&&\\
&\leq&
\mbox{$\displaystyle \liminf_{n \to \infty}
\frac{\log \left(4 \mathbb{P}\left(\left|S_{n} \right| > x\sqrt{ng(\log n)} \right)\right)}{g(\log n)}$ ~(by (5.1))}\\
&&\\
&=&
\mbox{$\displaystyle \liminf_{n \to \infty}
\frac{\log \mathbb{P}\left(\left|S_{n} \right| > x\sqrt{ng(\log n)} \right)}{g(\log n)}$}\\
&&\\
&\leq&
\mbox{$\displaystyle \limsup_{n \to \infty}
\frac{\log \mathbb{P}\left(\left|S_{n} \right| > x\sqrt{ng(\log n)} \right)}{g(\log n)}$}\\
&&\\
&\leq&
\mbox{$\displaystyle 0~~\mbox{for all}~ x > 0 ~(\mbox{since}~ \mathbb{P}(A) \leq 1~ \mbox{for any event} ~A)$.}\\
\end{array}
\end{equation}
Since $\lim_{c \to \infty} \mathbb{E}Y^{2}(c) = \mathbb{E}Y^{2} = \infty$, we see that (2.7) follows from (5.2).

For Case III, applying Theorem 2.1 (i.e., (2.4)), we have
\[
\begin{array}{lll}
0 & = &
\mbox{$\displaystyle
 - \left(\frac{x^{2}}{2\sigma^{2}} \wedge \frac{0}{2^{\rho}} \right) $}\\
 &&\\
 & \leq &
 \mbox{$\displaystyle
 \liminf_{n \to \infty} \frac{\log \mathbb{P}\left(\left|S_{n} \right| > x\sqrt{ng(\log n)} \right)}{g(\log n)} $
 ~~(by Theorem 2.1)}\\
 &&\\
& \leq &
\mbox{$\displaystyle
\limsup_{n \to \infty} \frac{\log \mathbb{P}\left(\left|S_{n} \right| > x\sqrt{ng(\log n)} \right)}{g(\log n)}$}\\
&&\\
& = &
\mbox{$\displaystyle
- \left(\frac{x^{2}}{2 \sigma^{2}} \wedge \frac{0}{2^{\rho}} \right) $~~(by Theorem 2.1)}\\
&&\\
& = & 0 ~~\mbox{for all}~ x > 0
\end{array}
\]
which also ensures (2.7).

Obviously, (2.7) implies (2.8).

We now establish the implication (2.8) $\Rightarrow$ (2.9) by contradiction. If (2.9) does not hold, then either
\[
\mathbb{E}X = 0 ~~\mbox{and}~~\mathbb{E}X^{2} = 0
\]
or
\[
\mathbb{E}X = 0, ~~\mathbb{E}X^{2} = \sigma^{2} \in (0, \infty), ~ \mbox{and}~~\lambda_{1} \neq \lambda_{2}.
\]
For the first case, for each $n \geq 2$, we have
\[
\mathbb{P}\left(\left|S_{n} \right| > x \sqrt{n g(\log n)} \right) = 0~~\mbox{for all}~ x > 0
\]
and hence,
\[
\lim_{n \to \infty} \frac{\log \mathbb{P}\left(\left|S_{n} \right| > x\sqrt{ng(\log n)} \right)}{g(\log n)}
= - \infty ~~\mbox{for all}~ x > 0
\]
which is contradictory to (2.8).

For the second case, since $\lambda_{1} \neq \lambda_{2}$ and $0 \leq \lambda_{1} \leq \lambda_{2} \leq \infty$,
we have $0 < \lambda_{2} \leq \infty$ and hence, by Theorem 2.1,
\[
\liminf_{n \to \infty} \frac{\log \mathbb{P}\left(\left|S_{n} \right| > x\sqrt{ng(\log n)} \right)}{g(\log n)}
\leq - \left(\frac{x^{2}}{2 \sigma^{2}} \wedge \frac{\lambda_{2}}{2^{\rho}} \right) < 0 ~~\mbox{for all}~ x > 0
\]
which is also contradictory to (2.8). ~$\Box$

\vskip 0.5cm




\vskip 0.5cm

{\bf References}

\begin{enumerate}
\item Bahadur, R. R., Zabell, S. L.: Large deviations of the sample mean in general vector spaces.
    Ann. Probab. {\bf 7}, 587-621 (1979).

\item Bingham, N. H., Goldie, C. M., Teugels, J. L.: Regular Variation. Encyclopedia of Mathematics
    and Its Applications {\bf 27}. Cambridge Univ. Press. (1987)

\item Bolthausen, E.: On the probability of large deviations in Banach spaces, Ann. Probab. {\bf 12},
    427-435 (1984).

\item Book, S. A.: Large deviations and applications, Encyclopedia of Statistical Sciences. {\bf 4}
    (S. Kotz, N. L. Johnson, and C. B. Read, eds.), John Wiley \& Sons, New York, 476-480 (1983).

\item Borovkov, A. A., Mogul'ski\u{i}, A. A.: Probabilities of large deviations in topological spaces.
    I, Sibirsk. Mat. Zh. {\bf 19}, no. 5, 988-1004 (1978) (Russian), translated in Siberian Math.
    J. {\bf 19}, no. 5, 697-709 (1979).

\item Chen, X.: Probabilities of moderate deviations for B-valued independent random vectors,
    Chinese J. Contemp. Math. {\bf 11}, 381-393 (1990).

\item Chen, X.: Moderate deviations of independent random vectors in a Banach space, Chinese J.
    Appl. Probab. Statist. {\bf 7}, no. 1, 24-32 (1991) (Chinese).

\item Chernoff, H.: A measure of asymptotic efficiency for tests of a hypothesis based on the sum
    of observations. Ann. Math. Statistics {\bf 23}, 493-507 (1952).

\item Comman, H.: Improvements of Plachky-Steinebach theorem. Theory Probab. Appl. {\bf 63}, 117-134 (2018).

\item Cram\'{e}r, H.: Sur un nouveau th\'{e}or\`{e}me-limite de la th\'eorie des probabilit\'{e}s.
    Actualit\'{e}s Sci. Indust. {\bf 736}, 5-23 (1938).

\item de Acosta, A.: Moderate deviations and associated Laplace approximations for sums of independent
    random vectors, Trans. Amer. Math. Soc. {\bf 329}, no. 1, 357-375 (1992).

\item Dembo, A., Zeitouni, O.: Large Deviations Techniques and Applications. Springer, Berlin-Heidelberg
     (2009).

\item Donsker, M. D., Varadhan, S. R. S.: Asymptotic evaluation of certain Markov process
    expectations for large time. III. Comm. Pure Appl. Math. {\bf 29}, 389-461 (1976).

\item Eichelsbacher, P., Lo\"{w}e, M.: Moderate deviations for i.i.d. random variables.
    ESAIM - Probab. Stat. {\bf 7}, 209-218 (2003).

\item Gantert, N.: A note on logarithmic tail asymptotics and mixing. Statist. Probab. Lett.
    {\bf 49}, 113-118 (2000).

\item Hu, Y. J., Nyrhinen, H.: Large deviations view points for heavy-tailed random walks.
    J. Theoret. Probab. {\bf 17}, 761-768 (2004).

\item Ledoux, M: Sur les d\'{e}viations mod\'{e}r\'{e}es des sommes de variables al\'eatoires vectorielles ind\'{e}pendantes
    de m\^{e}me loi [On moderate deviations of sums of i.i.d. vector random variables], Ann. Inst. H. Poincar\'{e} Probab.
    Statist. {\bf 28}, no. 2, 267-280 (1992) (French).

\item Ledoux, M., Talagrand, M.: Probability in Banach Spaces:  Isoperimetry and Processes.
    Springer-Verlag, Berlin (1991).

\item Li, D., Miao, Y.: A supplement to the laws of large numbers and the large deviations.
     Stochastics {\bf 93}, 1261-1280 (2021).

\item Li, D., Miao, Y., Stoica, G.: A general large deviation result for partial sums of i.i.d.
    super-heavy tailed random variables. Stat. Prob. Lett. {\bf 184}, Article 109371 (2022).

\item Li, D., Rosalsky: Precise lim sup behavior of probabilities of large deviations for sums of
    i.i.d. random variables, Int. J. Math. Math. Sci. {\bf 2004}, 3565-3576 (2004).

\item Li, D., Rosalsky, A., Al-Mutairi, D. K.: A large deviation principle for bootstrapped sample
    means. Proc. Amer. Math. Soc. {\bf 130}, 2133-2138 (2002).

\item Nakata, T.: Large deviations for super-heavy tailed random walks. Stat. Prob. Lett. {\bf 180},
    Article 109240 (2022).

\item Petrov, V. V.: Generalization of Cram\'{e}r's limit theorem, Uspehi Matem. Nauk (N.S.) {\bf 9},
    no. 4(62), 195-202 (1954) (Russian), translated in Select. Transl. in Math. Stat. and
    Probab. {\bf 6}, 1-8 (1966).

\item Petrov, V. V.: A generalization of a certain inequality of L\'{e}vy, Teor. Veroyatnost. i Primenen.
    {\bf 20}, 140-144 (1975) (Russian), translated in Theory Probab. Appl., {\bf 20}, 141-145 (1975).

\item Petrov, V. V.: Sums of Independent Random Variables. Springer-Verlag, New York (1975).

\item Petrov, V. V.: Limit Theorems of Probability Theory. Sequences of Independent Random Variables.
     Oxford Studies in Probability, {\bf 4}, The Clarendon Press, Oxford University Press, New York (1995).

\item Plachky, D.: On a theorem of G. L. Sievers. Ann. Math. Statist. {\bf 42}, 1442-1443 (1971).

\item Plachky, D., Steinebach, J.: A theorem about probabilities of large deviations with an application to queuing theory. Period. Math. Hungar. {\bf 6}, 343-345 (1975).

\item Saulis, L., Statulevi\u{c}ius, V. A.: Limit Theorems for Large Deviations, Mathematics and its
    Applications (Soviet Series), {\bf 73}, Kluwer Academic Publishers, Dordrecht (1991).

\item Sievers, G. L.: On the probability of large deviations and exact slopes. Ann. Math. Statist. {\bf 40}, 1908-1921 (1969).

\item Stoica, G.: Large gains in the St. Petersburg game. C. R. Math. Acad. Sci. Paris {\bf 346},
    no. 9-10, 563-566 (2008).

\item Stroock, D. W.: An Introduction to the Theory of Large Deviations. Springer, New York (1984).

\end{enumerate}

\end{document}